\documentclass[11pt]{article}

\usepackage{lineno,dsfont}
\usepackage{amsmath, latexsym,amssymb,amsopn,amscd,amsbsy,vmargin,bm,marginnote,cprotect}
\usepackage{mathrsfs}
\usepackage{subfig}
\usepackage{natbib}
\usepackage[usenames,dvipsnames]{color}
\usepackage[normalem]{ulem}
\usepackage{graphicx}
\usepackage{float}
\usepackage{epsfig}
\usepackage{relsize}
\usepackage{bigints}
\usepackage{xcolor}
\usepackage{bm}                   
\usepackage[latin1]{inputenc}  
\usepackage{amsthm} 

\newcommand{\N}{\ensuremath{\mathbb{N}}}
\newcommand{\R}{\ensuremath{\mathbb{R}}}

\newcommand{\Prob}{\ensuremath{\mathbb{P}}}
\newcommand{\Esp}{\ensuremath{\mathbb{E}}}
\newcommand{\A}{\ensuremath{\forall}}

\newcommand{\B}{\mathscr{B}}
\newcommand{\F}{\cl{F}}
\newcommand{\daw}{\downarrow}
\newcommand{\Ind}{\mathlarger{\mathbf{1}}}
\newcommand{\captionstring}[1]{\noexpand\noexpand\noexpand\string\string#1}

\newcommand{\cl}[1]{\mathcal{#1}}

\newcommand{\msf}[1]{\mathsf{#1}}
\newcommand{\mbf}[1]{\mathbf{#1}}

\def\l{\left}
\def\r{\right}

\def\deq{\stackrel{d}{=}}
\def\aseq{\stackrel{a.s.}{=}}
\def\dto{\stackrel{d}{\to}}
\def\asto{\stackrel{a.s.}{\to}}
\def\wto{\stackrel{w}{\to}}
\def\dwto{\stackrel{dw}{\to}}
\def\L2to{\stackrel{\cl{L}_2}{\to}}

\def\iid{\stackrel{\mbox{\scriptsize{iid}}}{\sim}}

\def\d{{\mbf{d}}}
\def\u{{\mbf{u}}}
\def\v{{\mbf{v}}}
\def\w{{\mbf{w}}}
\def\x{{\mbf{x}}}
\def\y{{\mbf{y}}}

\def\m{{\mbf{m}}}
\def\p{{\mbf{p}}}
\def\W{{\mbf{W}}}
\def\K{{\mbf{K}}}
\def\U{{\mbf{U}}}
\def\V{{\mbf{V}}}
\def\X{{\mbf{X}}}

\def\D{{\mbf{D}}}
\def\bmu{{\bm{\mu}}}
\def\bxi{{\bm{\xi}}}
\def\bXi{{\bm{\Xi}}}
\def\bkap{{\bm{\kappa}}}
\def\bsig{{\bm{\sigma}}}
\def\bphi{{\bm{\phi}}}
\def\bgamma{{\bm{\gamma}}}
\def\dw{{\mbf{w}}^{\daw}}
\def\tw{{\mbf{\tilde{w}}}}
\def\dW{{\mbf{W}}^{\daw}}
\def\tW{{\mbf{\tilde{W}}}}

\def\mathL{\cl{L}}

\newtheorem{theo}{Theorem}[section]
\newtheorem{prop}[theo]{Proposition}
\newtheorem{lem}[theo]{Lemma}
\newtheorem{rem}[theo]{Remark}
\newtheorem{cor}[theo]{Corollary}

\title{\bf Beta-Binomial stick-breaking non-parametric prior}

\author{
  Mar\'ia F. Gil--Leyva\\
  IIMAS, Universidad Nacional Aut\'onoma de M\'exico\\
  CDMX, M\'exico\\
  \texttt{marifer@sigma.iimas.unam.mx} \\
  \and
  Rams\'es H. Mena\\
  IIMAS, Universidad Nacional Aut\'onoma de M\'exico\\
  CDMX, M\'exico\\
  \texttt{ramses@sigma.iimas.unam.mx} \\
  \and
  Theodoros Nicoleris\\
  Department of Economics\\ 
  National and Kapodistrian University of Athens\\
  Athens, Greece\\
  \texttt{tnicoleris@econ.uoa.gr}
}

\date{}

\def\spacingset#1{\renewcommand{\baselinestretch}%
{#1}\small\normalsize} \spacingset{1}

\begin{document}
\maketitle

\begin{abstract}
A new class of nonparametric prior distributions, termed Beta-Binomial stick-breaking process, is proposed. By allowing the underlying length random variables to be dependent through a Beta marginals Markov chain, an appealing discrete random probability measure arises. The chain's dependence parameter controls the ordering of the stick-breaking weights, and thus tunes the model's label-switching ability. Also, by tuning this parameter, the resulting class contains the Dirichlet process and the Geometric process priors as particular cases, which is of interest for MCMC implementations.  

Some properties of the model are discussed and a density estimation algorithm is proposed and tested with simulated datasets.
\end{abstract}

\noindent%
\textbf{{\it Keywords:}} Beta-Binomial Markov chain, Density estimation, Dirichlet process prior, Geometric process prior, Stick-breaking prior
\vfill

\newpage
\spacingset{1.1}

\section{Introduction}
Discrete random probability measures and their distributions play a key role in Bayesian nonparametric statistics. The availability of general classes of priors and their different representations are crucial for the study of theoretical properties, as well as for 
the proposal of  simulation and estimation algorithms. This continuously encourages the search of competitive  alternatives to the canonical model, \cite{F73} Dirichlet process. At the outset, one could consider a (proper) species sampling process \citep[][]{P06} over a measurable Polish space $(S,\B(S))$, %which takes the form
\begin{equation}\label{eq:s.s.p}
\bmu = \sum_{j \geq 1}\w_j\delta_{\bxi_j},
\end{equation}where the atoms, $\bXi = \l(\bxi_j\r)_{j \geq 1}$, and the weights, $\W = \l(\w_j\r)_{j \geq 1}$, are independent collections of  random variables (r.v.'s), with $\bxi_j\iid P_0$, a diffuse measure on $(S,\B(S))$,  and $\sum_{j \geq 1}\w_j = 1$, almost surely (a.s.). To fully specify the law of $\bmu$, one could assume a form for $P_0$ and place a distribution over the infinite dimensional simplex $\Delta_{\infty} = \{(w_1,w_2,\ldots): w_i \geq 0, \sum_{i \geq 1}w_i = 1\}$. An important aspect to note is that
\begin{equation}\label{eq:label_switch}
\sum_{j \geq 1}\w_{j}\delta_{\bxi_j} \deq \sum_{j \geq 1}\w_{\rho(j)}\delta_{\bxi_j}
\end{equation}
for every random permutation of $\N$, $\rho$, independent of $\bXi$. This means that once the atom's distribution, $P_0$, is fixed, there are infinitely many distributions over $\Delta_\infty$ that lead to the exact same prior, hence the need to study orderings for the weights. In particular, one can consider the decreasing ordering of its elements, here denoted by $\dW = (\dw_j)_{j \geq 1}$, with $\dw_1 > \dw_2 > \cdots$ a.s.,  or the size-biased permutation, denoted by $\tW = (\tw_j)_{j \geq 1}$, which satisfies $\Prob[\tw_1 = \w_j|\W] = \w_j$, and for $n \geq 2$
\[
\Prob[\tw_n = \w_j|\W,\tw_1,\ldots\tw_{n-1}] =\frac{\w_j}{1-\sum_{i=1}^{n-1}\tw_{i}} \Ind_{\{\w_j \not\in \{\tw_1,\ldots,\tw_{n-1}\}\}}.
\]
Working with  decreasing representations of the weights reduces the identifiability problem that arises from \eqref{eq:label_switch} in the sense that if $\bgamma_1,\bgamma_2,\ldots$ is sampled i.i.d. from $\bmu$, conditionally given $\bmu$, then $\dw_1$ corresponds to the atom that appears more frequently in the sequence, $\dw_2$ corresponds to the second most frequent value, and so on \citep[e.g.,][]{MW15}. On the other hand, the size-biased permutation of the weights is of interest when the focus is in the clusters featured in the sample,  i.e. if $\bgamma^{*}_j$ is the $j$th distinct value to appear in the sample, then the long-run proportion of elements in $\{n: \bgamma_n = \bgamma^{*}_j\}$ coincides precisely with $\tw_j$ \citep[][]{P96b}.\\

Different  techniques to place distributions on $\Delta_{\infty}$ are available \citep[e.g.][]{F73,BM73,JLP09} and connections among such techniques are well known \citep[e.g.][]{IJ01,IZ02,HHMW10}. Perhaps one of the most practical constructions is enjoyed by the so-called stick-breaking process \citep[][]{McCloskey65,S94,IJ01} where the weights are decomposed as
\begin{equation}\label{eq:s.b.s}
\w_1 = \v_1, \quad \w_j = \v_j\prod_{i=1}^{j-1}(1-\v_i), \quad j \geq 2,
\end{equation} for some sequence taking values in $[0,1]$, $\V = \l(\v_i\r)_{i \geq 1}$, hereinafter referred to as length variables (l.v.'s). 
%This representation translates the problem of specifying a distribution over $\Delta_\infty$ (law of $\W$) into placing one over $[0,1]^{\infty}$ (law of $\V$) and further corroborating that through \eqref{eq:s.b.s}, $\sum_{j \geq 1}\w_j = 1$ a.s., is satisfied. 
The practical compromise inherent to \eqref{eq:s.b.s} is relatively little, as most practical classes of priors have a stick-breaking representation,  e.g. the Dirichlet process \citep{F73,S94}, its two-parameter generalization \citep{PPY92,PY92}, the  normalized inverse-Gaussian process \citep{FLP12} and the more general class of homogeneous normalized random measures with independent increments \citep{FLNNPT16}. In particular, the Dirichlet process is recovered when $\v_i\iid\msf{Be}(1,\theta)$, for some $\theta > 0$, and, as shown by \cite{P96}, the resulting weights coincide with the corresponding size-biased permutation of them, an ideal feature for clustering \citep[][]{P96b}. A  different stick-breaking prior is the Geometric process, introduced by \cite{FMW10}. For this case, the decreasing ordering of the weights takes the form
\[
\w_j = \bm{\lambda}(1-\bm{\lambda})^{j-1}, \quad j \geq 1,
\]
for some $\bm{\lambda} \sim \msf{Be}(\alpha,\theta)$, with $\alpha,\theta > 0$. Here the random variables $(\v_i)_{i \geq 1}$ are completely dependent, indeed identical, unlike for the Dirichlet process. As mentioned above, the ordering of the weights, or lack of it, is of high relevance when using Bayesian nonparametric priors for density estimation and/or clustering. The dependence on only one random variable makes the Geometric process an attractive choice from a numerical point of view, and also makes it quite simple to generalize to non-exchangeable settings \citep[][]{FMW09,MRW11,HMNW16}. Furthermore, as shown by \cite{BO14}, both the Dirichlet and the Geometric processes have full support.

We propose a new class of stick-breaking distributions over $\Delta_\infty$, featured by  dependent l.v.'s driven by a strictly stationary Beta Markov chain, thus leading to a novel family of random probability measures, the Beta-Binomial stick-breaking (BBSB) priors. The Beta Markov chain in question has a dependence parameter which modulates the ordering of the corresponding weights, allowing BBSB priors to enjoy a good trade-off between weights identifiability and  mixing. For extreme values of the dependence parameter, we find that the Dirichlet process and the Geometric process priors are particular cases of our model. Furthermore, using an extension of the aforementioned result by \cite{BO14}, we will see that BBSB priors also have full support.

The remaining part of the article is organized as follows: In Section~\ref{sec:Beta_Bin_chain} we present the construction of the  Markov chain with $\msf{Be}(\alpha,\theta)$ marginals. Inhere, we also analyse some special and limiting cases that will subsequently allow to recover the Dirichlet and Geometric processes. This Markov chain then assembles in Section~\ref{sec:Beta_Bin_prior} a sequence of l.v.'s, thus leading to Beta-Binomial stick-breaking priors.  In Section~\ref{sec:dens_est_scheme} we derive a sampling scheme for density estimation and, in Section~\ref{sec:illustrations} we test it in simulated data. The proofs of the main results are deferred to the Appendix.

\section{Beta-Binomial Markov chain}\label{sec:Beta_Bin_chain}

%:
Following \cite{PCW02}, given a density function $\pi_{\v,\x}(v,x)$ with marginals $\pi_{\v}(v)$ and $\pi_{\x}(x)$, and whose conditional distributions are $\pi_{\v|\x}(v|x)$ and $\pi_{\x|\v}(x|v)$, it is possible to construct two of reversible Markov chains $(\v_i)_{i \geq 1}$ and $(\x_i)_{i \geq 1}$ with stationary distributions $\pi_{\v}$ and $\pi_{\x}$ respectively. The construction considers the law induced by  $\v_1 \sim \pi_{\v}$, and  
$\{\x_i\mid \v_{i}\} \sim \pi_{\x|\v}(\cdot|\v_{i})$, $\{\v_{i+1} \mid \x_i\} \sim \pi_{\v|\x}(\cdot|\x_i)$, for $i \geq 1$; where $\v_{i+1}$ is conditionally independent of $(\v_1,\x_1,\ldots,\v_{i-1},\x_{i-1},\v_i)$ given $\x_i$, and analogously $\x_{i+1}$ is conditionally independent of $(\v_1,\x_1,\ldots,\v_i,\x_i)$ given $\v_{i+1}$. Arising from the Beta-Binomial conjugate model, we take
\[
\pi_{\v,\x}(v,x) = \msf{Bin}(x|\kappa,v)\msf{Be}(v|\alpha,\theta),
\]
for some $\alpha,\theta > 0$, $\kappa \in \{0,1,\ldots\}$, and where $\msf{Bin}(0,p) = \delta_0$. Thus, the dependence induced by $\v_1 \sim \msf{Be}(\alpha,\theta)$, and $\{\x_i\mid\v_{i}\} \sim \msf{Bin}(\kappa,\v_{i})$, $\{\v_{i+1}\mid\x_i\} \sim \msf{Be}(\alpha+\x_i,\theta+\kappa-\x_i)$, for $i \geq 1$ generates  Markov chains, $\V = (\v_i)_{i \geq 1}$ and $\X = (\x_{i})_{i \geq 1}$, where the former has transition probabilities given by
\begin{equation}\label{eq:v_trans}
\Prob[\v_i \in A|\v_{i-1}] = \int_A\sum_{x=0}^{\kappa} \msf{Be}(s|\alpha+x,\theta+\kappa-x)\msf{Bin}(x|\kappa,\v_{i-1})ds,
\end{equation}
and stationary distribution $\msf{Be}(\alpha,\theta)$, and the latter
\begin{equation}\label{eq:x_trans}
\begin{aligned}
\Prob[\x_i = x|\x_{i-1}] & = \int_0^1 \msf{Bin}(x|\kappa,p)\msf{Be}(p|\alpha+\x_{i-1},\theta+\kappa-\x_{i-1})dp\\
& = \binom{\kappa}{x}\frac{(\alpha+\x_{i-1})_{x\uparrow}(\theta+\kappa-\x_{i-1})_{\kappa-x\uparrow}}{(\alpha+\theta+\kappa)_{\kappa \uparrow}},
\end{aligned}
\end{equation}
where $(y)_{m \uparrow} = \prod_{j=0}^{m-1}(y+j)$, and its stationary distribution is
\begin{equation}\label{eq:x_station}
\Prob[\x_i = x] = \binom{\kappa}{x}\frac{(\alpha)_{x\uparrow}(\theta)_{\kappa-x\uparrow}}{(\alpha+\theta)_{\kappa \uparrow}}.
\end{equation}

To any Markov chains, $\V$, $\X$ and $(\V,\X) = (\v_i,\x_i)_{i \geq 1}$, we refer to them as 
\emph{Beta}, \emph{Binomial} and \emph{Beta-Binomial} chains. See \cite{NW02} and \cite{MW09} for more on this kind of Markov chains. In what follows, we focus on the the Beta chain and some of its properties, specifically in how the parameter $\kappa$ affects the dependence of the chain. This will be relevant for our construction of the nonparametric prior in the following section.

\begin{prop}\label{prop:prop}
Let $(\V,\X)$ be a Beta-Binomial chain with parameters $(\kappa,\alpha,\theta)$, then for the Beta chain, $\V$, and for every $i \geq 1$, we have the following conditional moments
\begin{itemize}
\item[\emph{a)}] 
$\displaystyle 
\Esp[\v_{i+1}|\v_i] = \frac{\alpha+\kappa\v_i}{\alpha+\theta+\kappa}.
$
\item[\emph{b)}] 
$\displaystyle
\msf{Var}(\v_{i+1}|\v_i) = \frac{(\alpha+\kappa\v_i)(\theta+\kappa(1-\v_i))+\kappa\v_i(1-\v_i)(\alpha+\theta+\kappa)}{(\alpha+\theta+\kappa)^2(\alpha+\theta+\kappa+1)}.
$
\item[\emph{c)}]
$\displaystyle
\msf{Cov}(\v_i,\v_{i+1}) = \frac{\kappa\alpha\theta}{(\alpha+\theta)^2(\alpha+\theta+1)(\alpha+\theta+\kappa)}.
$
\item[\emph{d)}]
$\displaystyle
\rho_{\v_i,\v_{i+1}} = \frac{\msf{Cov}(\v_i,\v_{i+1})}{\sqrt{\msf{Var}(\v_i)}\sqrt{\msf{Var}(\v_{i+1})}} = \frac{\kappa}{\alpha+\theta+\kappa}.
$
\end{itemize}
\end{prop}
Fixing the value of $\kappa$ and increasing either $\alpha$ or $\theta$, the correlation coefficient, $\rho_{\v_i,\v_{i+1}}$ goes to $0$. Conversely, if we fix $\alpha$ and $\theta$, for large values of $\kappa$, $\rho_{\v_i,\v_{i+1}} \approx 1$. Also,  if $\alpha$ and $\theta$ are very small with respect to $\kappa$
\[
\Esp[\v_{i+1}|\v_i] \approx \v_i\quad\mbox{and}\quad \msf{Var}(\v_{i+1}|\v_i) \approx \frac{2\v_i(1-\v_i)}{\kappa+1}.
\]

\begin{figure}[H]
\centering
\includegraphics[scale=0.4]{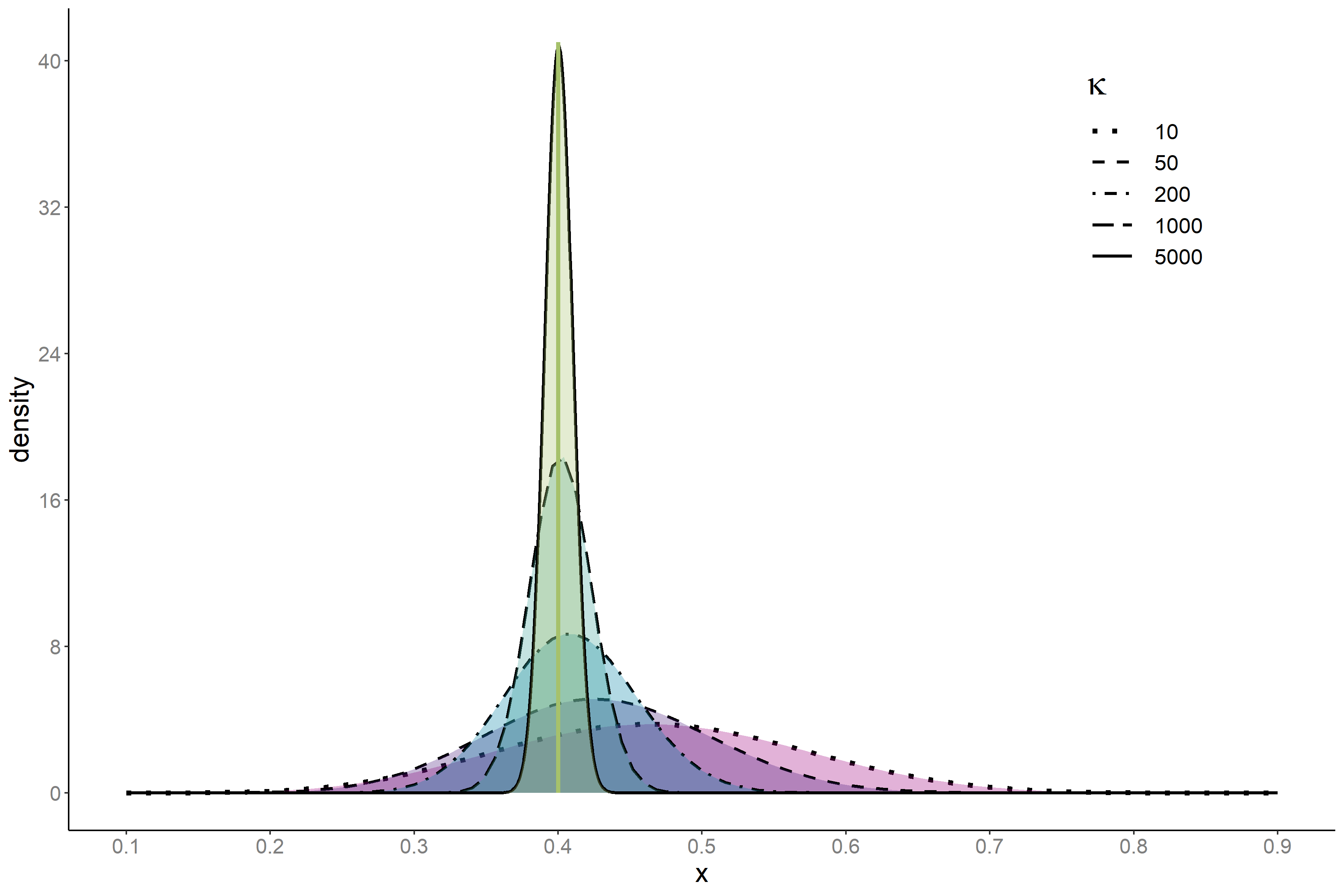} %I edit/change this file
\caption{Conditional densities of $\v_{i+1}$ given $\v_{i} = 0.4$, for distinct values of $\kappa$. We vary $\kappa$ in the set $\{10,50,200,1000,5000\}$,  the values of $\alpha$ and $\theta$ were fixed to $10$.\label{FigKappa}}
\end{figure}

Hence, intuition tells us that the conditional distribution of $\v_{i+1}$ given $\v_{i}$, tends to $\delta_{\v_{i}}$, as $\kappa$ grows, see  Figure~\ref{FigKappa}. The following result generalizes this intuition. 

\begin{prop}\label{prop:BB_DP_G}
Let $\V^{(\kappa)} = \l(\v_{i}^{(\kappa)}\r)_{i \geq 1}$ be a Beta-chain with parameters $(\kappa,\alpha,\theta)$.
\begin{itemize}
\item[\emph{(i)}] For $\kappa =0$, $\V^{(0)}$ is a sequence of i.i.d. random variables with distribution $\msf{Be}(\alpha,\theta)$.
\item[\emph{(ii)}] As $\kappa \to \infty$, $\V^{(\kappa)}$ converges in distribution to $(\bm{\lambda},\bm{\lambda},\ldots)$, where $\bm{\lambda} \sim \msf{Be}(\alpha,\theta)$.
\end{itemize}
\end{prop}

\section{Beta-binomial stick-breaking prior}\label{sec:Beta_Bin_prior}
We call Beta-Binomial stick-breaking prior to any species sampling process, $\bmu$,  with weights sequence as in \eqref{eq:s.b.s} for some l.v.'s, $\V$,  driven by a Beta chain with transition density \eqref{eq:v_trans}. As usual, the parameters of the l.v.'s are inherited to the prior, adding to the latter, the diffuse probability measure, $P_0$, as an additional parameter. The first property to check is that the corresponding weights add up to one. 
\begin{prop}\label{prop:BB_s.b.s.}
Let $\W$ be as in equation \eqref{eq:s.b.s}, for some Beta chain, $\V$. Then
\[
\sum_{j \geq 1} \w_j \aseq 1.
\]
\end{prop}
Moreover, notice that for every $0<\delta<\varepsilon<1$ and $n \geq 1$, any Beta-Binomial chain, $(\V,\X)$, with parameters $(\kappa,\alpha,\theta)$, satisfies
\begin{align*}
\Prob\l[\bigcap_{i=1}^n (\delta < \v_i < \varepsilon)\r] & = \Esp\l[\prod_{i=1}^n\Prob\l[\delta < \v_i < \varepsilon|\X\r]\r]\\
& = \Esp\l[\Prob[\delta < \v_1 < \varepsilon|\x_1]\prod_{i=2}^n\Prob\l[\delta < \v_i < \varepsilon|\x_{i-1},\x_i\r]\r] > 0,
\end{align*}
as conditionally given $\X$, the elements of $\V$ are independent and Beta distributed. As shown by \cite{BO14}, the above observation shows that any Beta-Binomial prior has full support, and thus feasible for nonparametric inference. The following results, which follow from Proposition~\ref{prop:BB_DP_G}, motivate their study.
\begin{theo}\label{theo:BB_DP_G}
Let $\bmu^{(\kappa)}$ be a BBSB prior  with parameters $(\kappa,\alpha,\theta,P_0)$ then
\begin{itemize}
\item[\emph{(i)}] For $\kappa = 0$ and $\alpha = 1$, $\bmu^{(0)}$ is a Dirichlet process with parameters $(\theta,P_0)$.
\item[\emph{(ii)}] For any $\alpha$ and $\theta$ fixed, as $\kappa \to \infty$, $\bmu^{(\kappa)}$ converges in distribution to the Geometric process, $\bmu$, with parameters $(\alpha,\theta,P_0)$.
\end{itemize}
\end{theo}

In terms of the ordering of the corresponding weights, we have the following corollary.

\begin{cor}\label{cor:dec_prob_weights}
Let $\l(\w_j^{(\kappa)}\r)_{j \geq 1}$ be as in equation \eqref{eq:s.b.s}, for some Beta chain, $\l(\v_i^{(\kappa)}\r)_{i \geq 1}$, with parameters $(\kappa,\alpha,\theta)$. Then
\begin{itemize}
\item[\emph{(i)}] For $\alpha = 1$, $\kappa = 0$, and any choice of $\theta$, $\l(\w_j^{(\kappa)}\r)_{j \geq 1}$ is size-biased ordered.
\item[\emph{(ii)}] For any choices of $\alpha$ and $\theta$, and for every $j \geq 1$
\[
\lim_{\kappa \to \infty}\Prob\l[\w_{j+1}^{(\kappa)} < \w_{j}^{(\kappa)}\r] = 1.
\]
\end{itemize}
\end{cor}

\begin{figure}[H]
\centering
\includegraphics[scale=0.4]{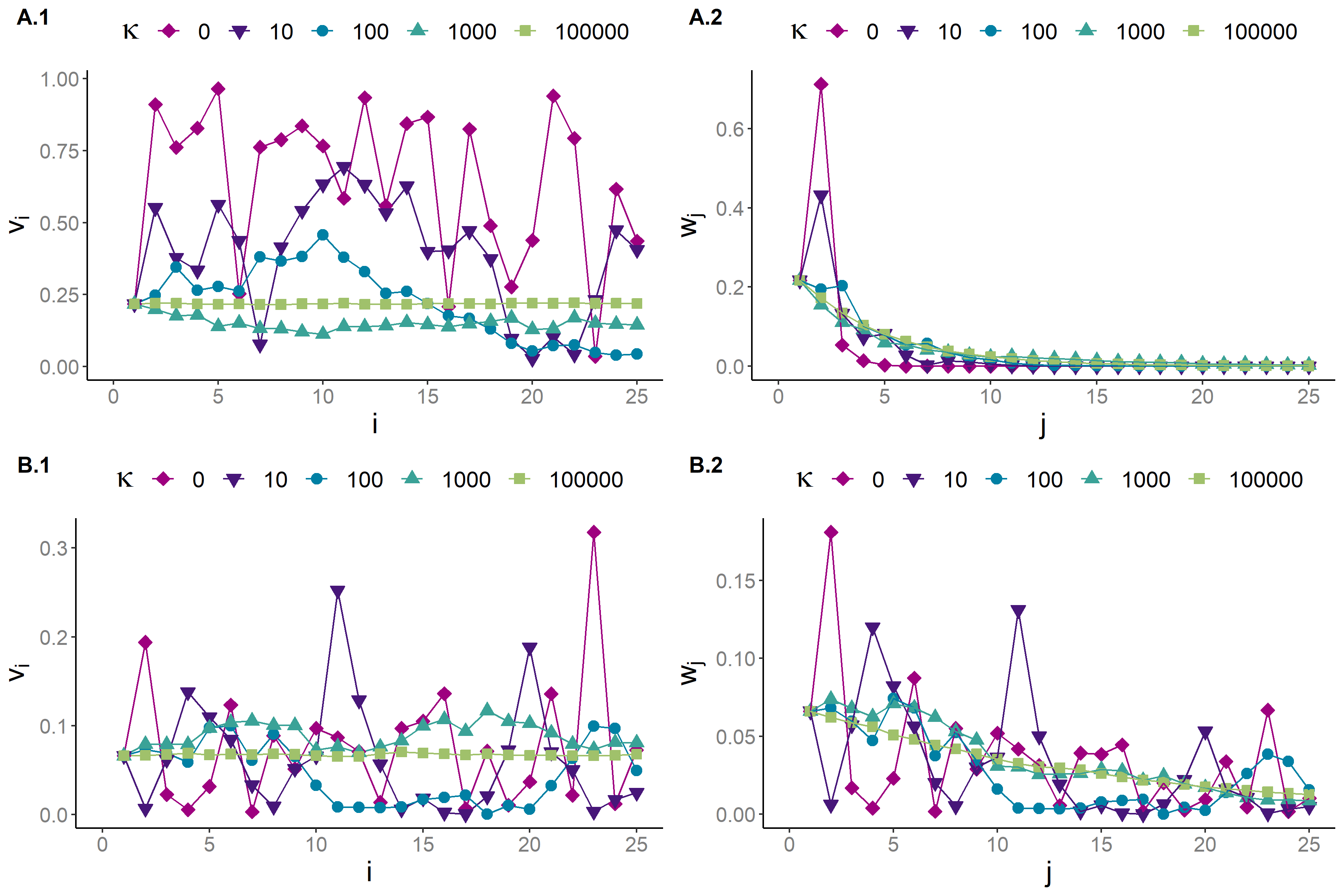}
\caption{\label{fig:weights} Simulations of $(\w_j)_{j=1}^{25}$ ($\msf{A.2}$ and $\msf{B.2}$) and their corresponding l.v.'s ($\msf{A.1}$ and $\msf{B.1}$ respectively) for distinct values of $\kappa$. For the Beta chains in $\msf{A.1}$, we fixed $\alpha = 1$ and $\theta = 1$, for the ones in $\msf{B.1}$ we used the same value of $\alpha$, whilst $\theta = 10$. The chains in a single graph share the same initial r.v. for the sake of a simpler analysis.}
\end{figure}

If we fix $\alpha = 1$, the choice $\kappa = 0$ implies that $\W = \tW$ is size-biased ordered. In general for such sequences $\Esp[\tw_j] \geq \Esp[\tw_{j+1}]$, even though $\tw_j \geq \tw_{j+1}$ does not occur with probability $1$. On the other extreme, as $\kappa \to \infty$ we have the decreasing ordering of the Geometric weights $\W = \dW$, which satisfy $\Prob\l[\dw_j \geq \dw_{j+1}\r] = 1$. Roughly speaking, by increasing the parameter $\kappa$, we make the weights sequence more likely to be decreasingly ordered. Figure \ref{fig:weights} shows some simulations of $(\w_j)_{j=1}^{25}$  and their corresponding l.v.'s  that illustrate the aforementioned behaviour.  The initial value, $\v_1$, of the Beta chain strongly affects  the behaviour of the complete sequence of weights, this is particularly evident for large values of $\kappa$. Recall that if $\kappa$ is sufficiently large we have $\v_1 \approx \v_2$, so for instance if $\v_1$ is close to $0$, then $(1-\v_1) \approx 1$ and $\w_2 = \v_2(1-\v_1) \approx \v_1 = \w_1$, which means that if $\v_2 > \v_1$ even slightly, we might obtain $\w_2 > \w_1$. Alternatively, a large value of $\v_1$, translates to a small value of $(1-\v_1)$, so in order to obtain $\w_2 > \w_1$, it would require $\v_2$ to be significantly larger than $\v_1$, which under the assumption that $\kappa$ is large, is not very likely to happen, as $\v_1 \approx \v_2$. The same intuition is inherited to the subsequent indexes since we also have $\v_2 \approx \v_3\approx \cdots$, for large values of $\kappa$. Hence,  the larger/smaller $\v_1$ is, the larger/smaller we expect $\v_i$ to be, for $i > 1$. Moreover, for large values of the parameter $\theta$ we expect $\v_1$ to take small values, thus in general,  a bigger value of $\theta$ requires an even larger value of $\kappa$, to induce a stochastically decreasing ordering of the weights.

\subsection{Distribution of the number of groups}

When working with any species sampling process, $\bmu$, such as a Dirichlet, BBSB or Geometric process\ldots, a r.v. of interest is the number of distinct values, $\K_n$, that a sample $\{\bgamma_1,\ldots,\bgamma_n\}$ driven by $\bmu$ exhibits. Although for some priors it is possible to compute or characterize the probabilistic behaviour of $\K_n$ \citep[see for instance][]{P06}, in general this is not an easy task to do. Despite this, whenever it is feasible to obtain samples from the weights sequence, $\W$, as is the case of any BBSB prior, obtaining samples from $\K_n$ can be easily achieved as follows: Sample $n$ independent $\msf{U}(0,1)$ r.v.'s, $(\u_k)_{k=1}^{n}$, and $(\w_j)_{j = 1}^{\varphi}$ where $\varphi$ is some constant satisfying $\sum_{j = 1}^{\varphi} \w_j > \max_k \u_k$. For $k \in \{1,\ldots,n\}$ and $i \in \{1,\ldots,\varphi\}$, let $\d_k = i$ if and only if $\sum_{j = 1}^{i-1} \w_j < \u_k < \sum_{j = 1}^{i} \w_j$ (with the convention that the empty sum equals $0$) then the number of distinct values $(\d_1,\ldots,\d_n)$ exhibits is precisely a sample from $\K_n$.

\begin{figure}[H]
\centering
\includegraphics[scale=0.4]{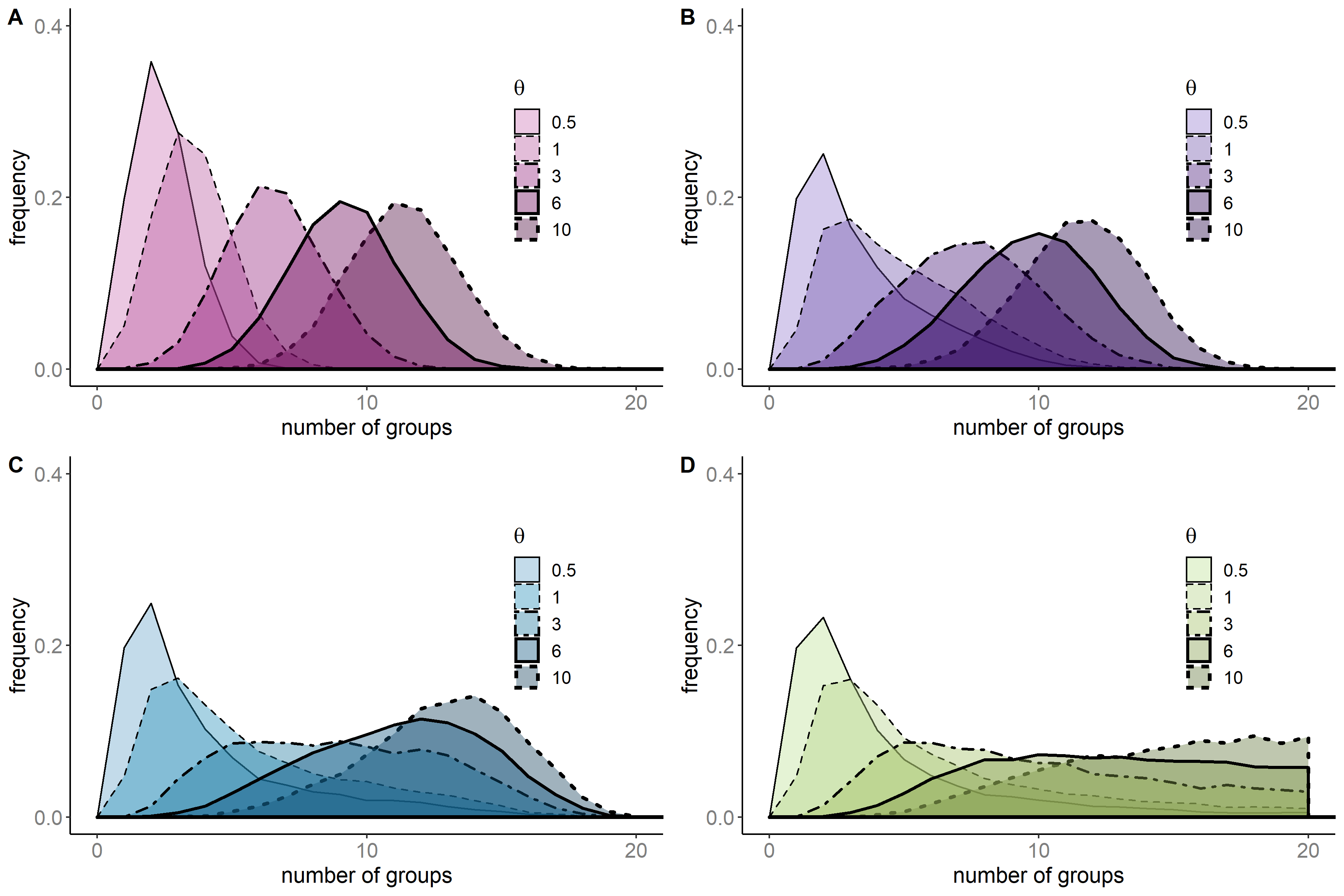}
\caption{\label{fig:Kna} Frequency polygons of samples of size $10000$ from $\K_{20}$ for distinct values of $\kappa$ and $\theta$ and fixing $\alpha = 1$. For the frequency polygons in $\msf{A},\msf{B}$ and $\msf{C}$ we fixed $\kappa$ to $0,10$ and $100$ respectively, whilst the frequency polygons in $\msf{D}$ correspond to the Geometric prior. For each fixed value of $\kappa$, we vary $\theta$ in the set $\{0.5,1,3,6,10\}$.}
\end{figure}

\begin{figure}[H]
\centering
\includegraphics[scale=0.4]{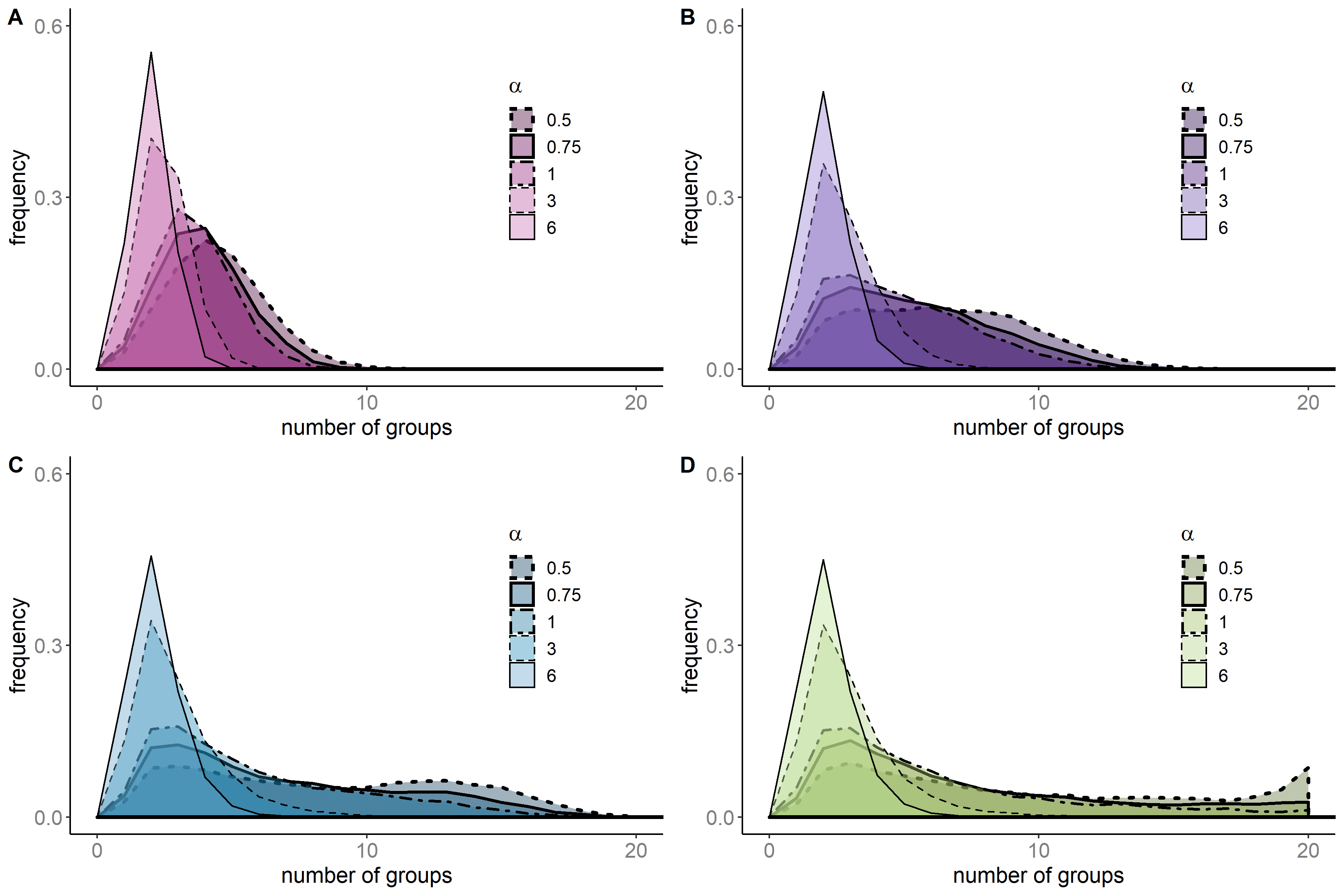}
\caption{\label{fig:Knb} Frequency polygons of samples of size $10000$ from $\K_{20}$ for distinct values of $\kappa$ and $\alpha$ and fixing $\theta = 1$. For the frequency polygons in $\msf{A},\msf{B}$ and $\msf{C}$ we fixed $\kappa$ to $0,10$ and $100$ respectively, whilst the frequency polygons in $\msf{D}$ correspond to the Geometric prior. For each fixed value of $\kappa$, we vary $\alpha$ in the set $\{0.5,0.75,1,3,6\}$.}
\end{figure}

To understand how the parameters of a BBSB prior affect the distribution of $\K_n$, we sampled as aforementioned varying the values of $\kappa$, $\alpha$ and $\theta$. Particularly, Figure \ref{fig:Kna}$(\msf{A})$ shows the distribution of $\K_n$ corresponding to the Dirichlet process, for which is well known that $\Esp[\K_n]$ increases when $\theta$ grows. This location behavior is also observed for other fixed values of $\kappa$ $(\msf{B},\msf{C}$ and $\msf{D})$. Figures \ref{fig:Kna} and \ref{fig:Knb}, illustrate how for fixed $\alpha$ and $\theta$, an increment on $\kappa$ contributes to the distribution of $\K_n$ with a heavier right tail, and thus a larger mean and variance, say the prior on $\K_n$ is less informative. In Figure \ref{fig:Kna}, where we fixed $\alpha = 1$, it can be observed that for bigger values of $\theta$, the distribution of $\K_n$ is more sensitive to an increment of $\kappa$. The same can be seen in Figure \ref{fig:Knb}, for fixed $\theta = 1$ and smaller values of $\alpha$.

\section{Density estimation for Beta-Binomial mixtures}\label{sec:dens_est_scheme}
Given a BBSB prior, $\bmu$, and a density kernel, $g(\cdot|s)$, with parameter space $S$, we can consider BBSB mixtures. Namely, we can model elements in $\y^{(n)} = \{\y_1,\ldots,\y_n\}$ as i.i.d. sampled from the random density

\begin{equation}\label{eq:f(y)_int}
\bphi(\y) := \pi(\y|\W,\bXi) = \int_{S}g(\y|s)\bmu(ds)=\sum_{j \geq 1}\w_jg(\y|\bxi_j).
\end{equation}

For MCMC implementation purposes, and following \cite{W07}, this random density can be augmented as
\begin{equation}\label{eq:f(y,u)}
\pi(\y,\u|\W,\bXi) = \sum_{j \geq 1} \Ind_{\{\u < \w_j\}}g(\y|\bxi_j),
\end{equation} 
where it can be easily deduced

\begin{equation}\label{eq:f(u|w)}
\pi(\u|\W) = \sum_{j \geq 1}\Ind_{\{\u<\w_j\}}.
\end{equation}
As in the Dirichlet process case, given $\u$, the number of components in the mixture is finite, with indexes being the elements of $A_\u(\W) = \{j: \u < \w_j\}$, that is 
\begin{equation}\label{eq:f(y|u,w,s)}
\pi(\y|\u,\W,\bXi) = \frac{1}{|A_\u(\W)|} \sum_{j \in A_\u(\W)} g(\y|\bxi_j).
\end{equation}
Using the membership variable $\d$, i.e. $\d = j$ iff $\y$ is sampled from $g(\cdot|\bxi_j)$, one can further consider the augmented joint density

\begin{equation}\label{eq:f(y,u,d)}
\pi(\y,\u,\d|\W,\bXi) = \Ind_{\{\u<\w_d\}}g(\y|\bxi_d).
\end{equation}
The complete data likelihood based on a sample of size $n$ from \eqref{eq:f(y,u,d)} is easily seen to be
\begin{equation}\label{eq:likelihood}
\mathL_{\bxi,\w}((\y_k,\u_k,\d_k)_{k=1}^{n}) = \prod_{k=1}^{n} \Ind_{\{\u_k < \w_{d_k}\}}g(\y_k|\bxi_{\d_k}),
\end{equation}
and under the assumption $P_0$ has a density, $p_0$, with respect to a suitable measure, the full joint density of every variable involved is
\begin{equation}\label{eq:joint}
\begin{aligned}
\pi(&(\y_k,\u_k,\d_k)_{k=1}^{n},(\v_i)_{i \geq 1},(\bxi_j)_{j \geq 1})\\
& = \left(\prod_{k=1}^{n} \Ind_{\{\u_k < \w_{\d_k}\}}g(\y_k|\bxi_{\d_k})\right) \prod_{j \geq 1}p_0(\bxi_j)\\
& \quad \quad \quad \times \left(\msf{Be}(\v_1|\alpha,\theta)\prod_{i \geq 1}\sum_{x=0}^{\kappa}\msf{Be}(\v_{i+1}|\alpha+x,\theta+\kappa-x)\msf{Bin}(x|\kappa,\v_i) \right),
\end{aligned}
\end{equation}
recall $\w_{d_k} = \v_{d_k}\prod_{i = 1}^{\d_k-1}(1-\v_i)$ with the convention that the empty product equals $1$.
\subsection{Full conditionals}
The full conditional distributions, required for posterior inference via a Gibbs sampler implementation, are proportional  to \eqref{eq:joint}, and given as follows.

\subsubsection*{1. Updating $\bXi$:}
\begin{equation*}\label{eq:upd_theta_j}
\pi(\bxi_j|\ldots) \propto p_0(\bxi_j)\prod_{k \in D_j}g(\y_k|\bxi_j), \quad j \geq 1,
\end{equation*}
where $D_j = \{k \geq 1: \d_k = j\}$. If $p_0$ and $g$ form a conjugate pair, the above is easy to sample from.
\subsubsection*{2. Updating $\V$ and $\U = (\u_k)_{k=1}^{n}$ as a block:}
\begin{align*}
\pi(\U,\V|\ldots) & \propto \left(\prod_{k=1}^{n} \w_{\d_k}^{-1}\Ind_{\{\u_k<w_{\d_k}\}}\w_{\d_k}\right)\times\\
& \quad \quad \times\left(\msf{Be}(\v_1|\alpha,\theta)\prod_{i \geq 1} \sum_{x=0}^{\kappa}\msf{Be}(\v_{i+1}|\alpha+x,\theta+\kappa-x)\msf{Bin}(x|\kappa,\v_i) \right).
\end{align*}
As $\w_{\d_k} = \v_{\d_k}\prod_{i=1}^{\d_k-1}(1-\v_i)$, with the convention $\prod_{i=1}^{0} (\cdot) = 1$, then
\begin{align*}
\pi(\U,&\V|\ldots) \propto \left[\prod_{k=1}^{n} \w_{\d_k}^{-1}\Ind_{\{\u_k<\w_{\d_k}\}}\right] \left[\v_1^{\alpha_1}(1-\v_1)^{\theta_1}\msf{Be}(\v_1|\alpha,\theta)\right]\times\\
& \times\left[\prod_{i \geq 1}\sum_{x=0}^{\kappa} (\v_{i+1})^{\alpha_{i+1}}(1-\v_{i+1})^{\theta_{i+1}}\msf{Be}(\v_{i+1}|\alpha+x,\theta+\kappa-x) \msf{Bin}(x|\kappa,\v_i)\right]\\
\end{align*}
where 
\[
\alpha_i = \sum_{k=1}^n\Ind_{\{\d_k = i\}} \quad \text{ and } \quad \theta_i = \sum_{k=1}^n\Ind_{\{\d_k > i\}}.
\]
Recalling that for $m \in \N$, and $z>0$, $\Gamma(m+z) = (z)_{m\uparrow}\Gamma(z)$, we obtain
\begin{align*}
\pi&(\U,\V|\ldots) \propto \left[\prod_{k=1}^{n} \msf{U}(\u_k|0,\w_{\d_k})\right] \left[\msf{Be}(\v_1|\alpha+\alpha_1,\theta+\theta_1)\right]\\
& \times\left[\prod_{i \geq 1}\sum_{x=0}^{\kappa} \msf{Be}(\v_{i+1}|\alpha_{i+1}+\alpha+x,\theta_{i+1}+\theta+\kappa-x) \right.\\
& \times\left.\frac{(\alpha+x)_{\alpha_{i+1}\uparrow}(\theta+\kappa-x)_{\theta_{i+1}\uparrow}}{(\alpha+\theta+\kappa)_{(\alpha_{i+1}+\theta_{i+1})\uparrow}}\msf{Bin}(x|\kappa,\v_i)\right],\\
\end{align*}
with the convention $(z)_{0\uparrow} = 1$. Thus, to update $\V$ and $\U$, we first sample $\V$ from
\begin{align*}
\pi&(\V|\ldots(\text{exclude } \U)\ldots) \propto \left[\msf{Be}(\v_1|\alpha+\alpha_1,\theta+\theta_1)\right]\\
& \times\left[\prod_{i \geq 1}\sum_{x=0}^{\kappa} \msf{Be}(\v_{i+1}|\alpha_{i+1}+\alpha+x,\theta_{i+1}+\theta+\kappa-x) \right.\\
&\left.\times\frac{(\alpha+x)_{\alpha_{i+1}\uparrow}(\theta+\kappa-x)_{\theta_{i+1}\uparrow}}{(\alpha+\theta+\kappa)_{(\alpha_{i+1}+\theta_{i+1})\uparrow}}\msf{Bin}(x|\kappa,\v_i)\right],\\
\end{align*}
which can be normalized to a product of Beta densities mixtures, and latter sample $\U$ from $\pi(\U|\ldots) = \prod_{k=1}^{n} \msf{U}(\u_k|0,\w_{\d_k})$.

\subsubsection*{3. Updating $\D = (\d_k)_{k=1}^{n}$:}
\begin{equation*}\label{eq:upd_d_k}
\pi(\d_k = j|\ldots) \propto g(\y_k|\bxi_j)\Ind_{\{\u_k < \w_j\}}, \quad k \in \{1,\ldots,n\},
\end{equation*}
which is a discrete distribution with finite support, hence easy to sample from.

\begin{rem}[For the updating of  $\bXi$ and $\V$]
As it is well-known for this algorithm, we do not need to sample $\v_j$ and $\bxi_j$ for every $j \geq 1$, it suffices to sample enough of them so that step \emph{3} can take place. Explicitly, it suffices to sample $\bxi_j$ and $\v_j$ for $j \leq \varphi$,  where $\varphi$ is a constant such that $\sum_{j=1}^{\varphi} \w_j \geq \max_k(1-\u_k)$, then it is not possible that $\w_j > \u_k$ for any $k \leq n$ and $j > \varphi$. 
\end{rem}

\subsection{Posterior distribution analysis}

Given $\l\{\l(\bxi^{(t)}_j\r)_j,\l(\w^{(t)}_j\r)_j,\l(\u^{(t)}_k\r)_k,\l(\d^{(t)}_k\r)_k\r\}_{t = 1}^{T},$ from $\{\bXi,\W,\U,\D|\y^{(n)}\}$ obtained after $T$ iterations of the Gibbs sampler, following \eqref{eq:f(y|u,w,s)} we estimate the density of the data by
\begin{equation}\label{eq:hat_f}
\Esp\l[\bphi\big|\y^{(n)}\r] \approx \frac{1}{T}\sum_{t=1}^{T} \frac{1}{n}\sum_{k=1}^{n} \frac{1}{\big|A_k^{(t)}\big|} \sum_{j \in A_k^{(t)}} g\l(\cdot\big|\bxi^{(t)}_j\r),
\end{equation}
where $A_k^{(t)} = \l\{j: \u_k^{(t)} < \w_j^{(t)}\r\}$. Furthermore, we can also estimate the posterior distribution of $\{\K_n|\y^{(n)}\}$ through
\begin{equation}\label{eq:post_Kn}
\Prob\l[\K_n = m\big|\y^{(n)}\r] \approx \frac{1}{T} \sum_{t=1}^{T} \Ind_{\{\K_n^{(t)} = m\}},
\end{equation}
where $\K_n^{(t)}$ is the number of distinct values $\l(\d_k^{(t)}\r)_k$ exhibits. As usual, when working with mixtures of densities, $\K_n$ can be interpreted as the number of components of the mixture featured in the sample $\y^{(n)}$, that is the number of elements in $\{g(\cdot|\bxi_j)\}_{j \geq 1}$ such that $\y_k$ is sampled from $g(\cdot|\bxi_j)$, for some $\y_k \in \y^{(n)}$. This way, the estimates \eqref{eq:hat_f} together with \eqref{eq:post_Kn}, give us information of how well a model performs for the given data set. Among the models for which \eqref{eq:hat_f} adjusts well to the data, those for which \eqref{eq:post_Kn} favours smaller values of $m$ might be preferred, as this means the model is mixing the components, $\{g(\cdot|\bxi_j)\}_{j \geq 1}$, more efficiently.

\subsection{Posterior inference for the dependence parameter}

In order to highlight the role of the dependence parameter, $\kappa$, we incorporate its posterior inference. Namely, we consider this parameter random and endow it with a prior distribution, $\bkap\sim\pi_{\kappa}$. For this case, the likelihood \eqref{eq:likelihood} remains identical and the joint distribution  \eqref{eq:joint} is multiplied by  $\pi_{\kappa}(\bkap)$. It can easily be seen that, conditionally given $\bkap$, the full conditionals $\{\pi(\bxi_j|\ldots)\}_{j \geq 1}$,  $\{\pi(\d_k|\ldots)\}_{k=1}^n$  and $\pi(\V,\U|\ldots)$ also remain the same.  As to the full conditional of $\bkap$ given the rest of the r.v.'s, we have that
\begin{equation}\label{eq:post_kap}
\pi(\bkap = \kappa|\ldots) \propto \pi_{\kappa}(\kappa)\prod_{i \geq 1}\sum_{x=0}^{\kappa} \msf{Be}(\v_{i+1}|\alpha+x,\theta+\kappa-x) \msf{Bin}(x|\kappa,\v_i),
\end{equation} which is easy to sample from if $\pi_{\kappa}$ has finite support. Summarizing, at each iteration of the Gibbs sampler, we update $\bXi$, $\V$, $\U$ and $\D$ as above and add a fourth step in which we sample $\bkap$ from \eqref{eq:post_kap}.

Finally, given the samples $\l(\bkap^{(t)}\r)_{t=1}^{T}$ obtained after $T$ iterations of the Gibbs sampler, once the burn-in period has elapsed, we estimate the posterior distribution of $\bkap$ by
\[
\Prob[\bkap = \kappa|\y^{(n)}] \approx \frac{1}{T}\sum_{t=1}^T \Ind_{\{\bkap^{(t)} = \kappa\}}
\]

\section{Illustrations}\label{sec:illustrations}

In principle, every choice of $\kappa$ leads to robust posterior MCMC estimates, after an appropriate burn-in period and enough valid iterations. However, depending on the sample, initial conditions, and current parameter values in the Gibbs sampler, the need to more/less ordered weights, thus different values of $\kappa$, might be required. To test the performance of BBSB priors for density estimation, we first conduct a small experiment in which we fix the value of $\kappa$ to $0, 10, 100$ and $\infty$ and compare the results provided by the $4$ distinct models. Secondly, in order to choose the optimal value of $\kappa$ for a dataset and given that the rest of the parameters are fixed, we place a prior distribution on the dependence parameter and analyse its posterior distribution. Here we also compare our models to another well-known stick-breaking prior, the Pitman-Yor process \citep{PPY92,PY92}. In all cases we assume a Gaussian kernel with random location and scale parameters, i.e., for each $j \geq 1$, $\bxi_j = (\m_j,\p_j)$, and $g(\y|\bxi_j) = \msf{N}(\y|\m_j, \p_j^{-1})$. To attain a conjugate pair for $p_0$ and $g$, we assume $p_0(\bxi_j) = \msf{N}(\m_j|\vartheta,\tau\p_j^{-1})\msf{Ga}(\p_j|a,b)$, where $a = b = 0.5$, $\tau = 100$ and $\vartheta = n^{-1}\sum_{k = 1}^n \y_k.$

\subsection{Analysis for BBSB mixtures with fixed dependence parameter}

For this exercise we simulated a data set ($\msf{database}$ $\msf{1}$) containing $200$ observations and featuring 11 modes equally spaced.  As 
it is well known for this type of data, and if the parameter $\theta$ is not carefully chosen, the Dirichlet mixture under estimates the number modes featured in the sample.  Alternatively, Geometric mixtures do recognize every mode, but they tend to use a large number of mixture components. In order to study how BBSB priors perform in this context, and to compare them with the Dirichlet and Geometric processes, we fixed $\alpha = 1$, $\theta = 1$ and vary $\kappa$ in the set $\{0,10,100,\infty\}$. No burn-in period was considered, so that one may analyse the number of iterations required by the model to provide a good estimate. 

\begin{figure}[H]
\centering
\includegraphics[scale=0.3]{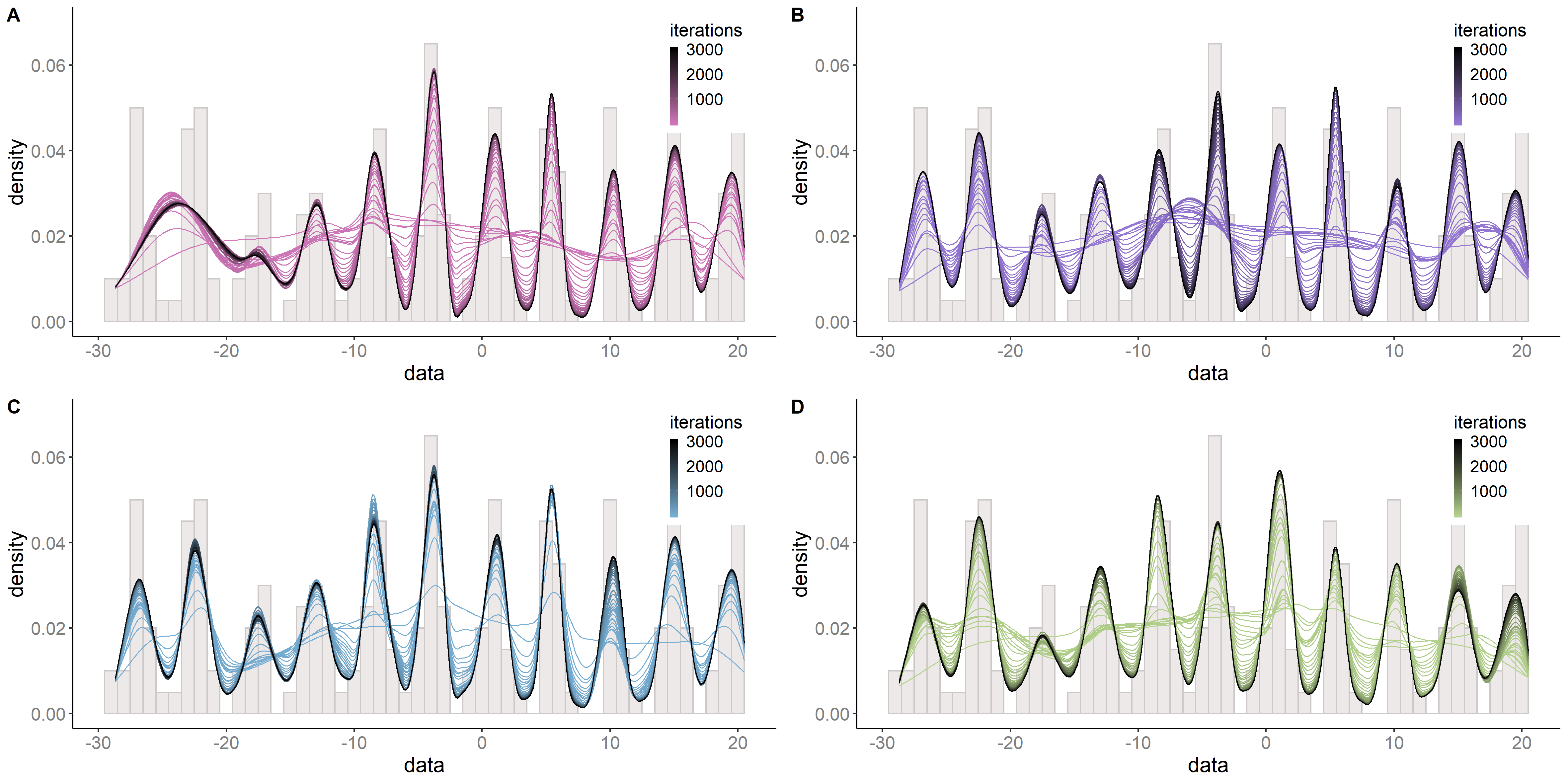}
\caption{\label{fig:dta_1} Evolution of the estimated densities for $\msf{database}$ $\msf{1}$, through the first $3000$ iterations of the Gibbs sampler, for four distinct BBSB mixtures. The estimated densities in $\msf{A},\msf{B},\msf{C}$ and $\msf{D}$ correspond to BBSB mixtures with $\kappa$ fixed to $0,10,100$ and $\infty$ respectively, in the four cases $\alpha = \theta = 1$.}
\end{figure}

\begin{figure}[H]
\centering
\includegraphics[scale=0.39]{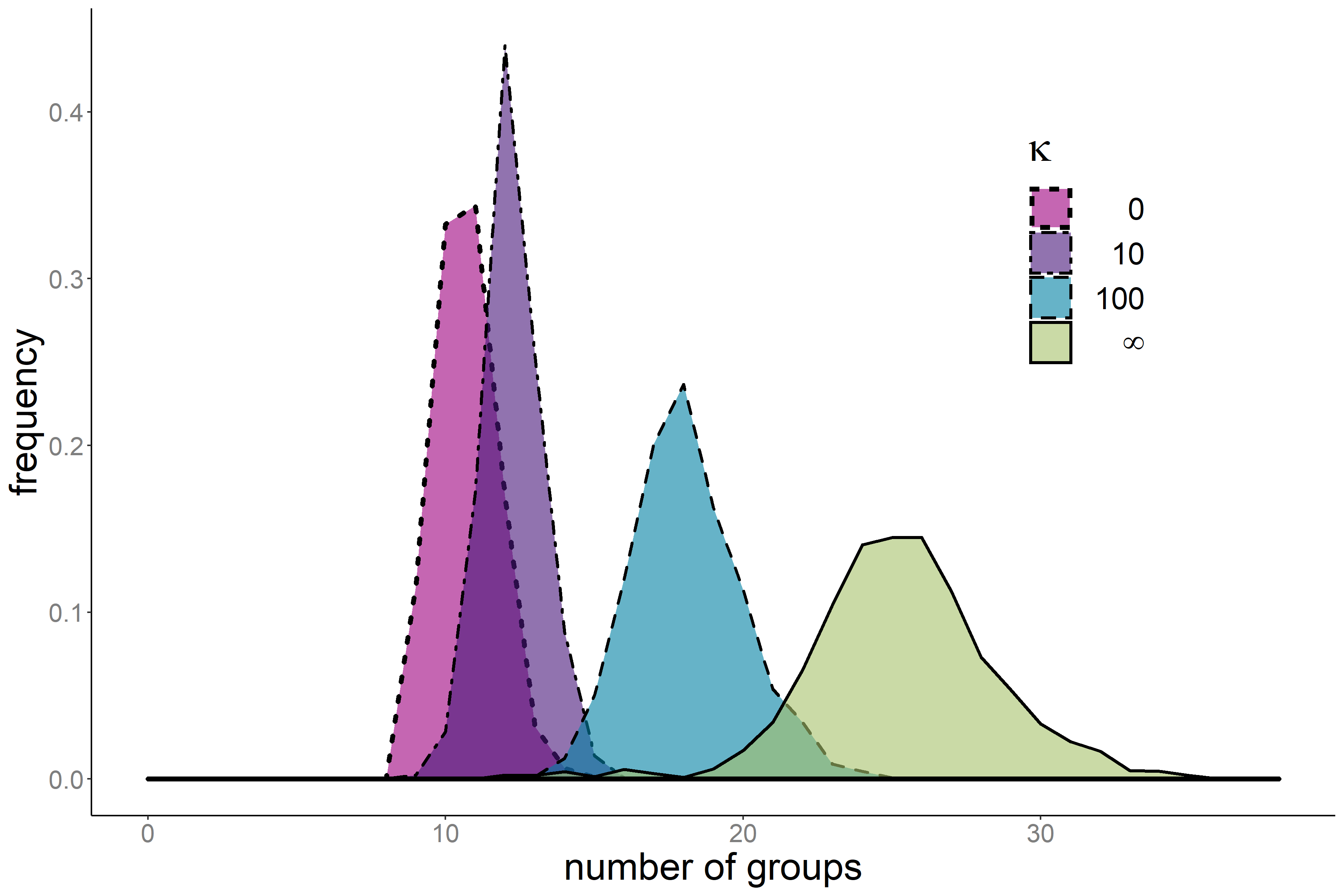}
\caption{\label{fig:post_Km_11} Frequency polygon of the estimated posterior distributions of $\K_n$ given $\msf{database}$ $\msf{1}$ for the four BBSB mixtures which share the parameters $\alpha = \theta = 1$, and differ on the parameter $\kappa$, same one that varies in the set $\{0,10,100,\infty\}$.}
\end{figure}

In Figure~\ref{fig:dta_1} we observe that the Dirichlet process ($\msf{A}$) fails to recover the eleven modes featured in the dataset, the three remaining models are able to capture the 11 well-separated modes. In terms of the number of iterations required to recognize the modes, we observe that BBSB mixtures with larger values of $\kappa$ ($\msf{C}$ and $\msf{D}$) perform better. Consistently with the prior analysis of the number of groups, in Figure~\ref{fig:post_Km_11} we observe that the posterior mean and variance increase as $\kappa$ does. Comparing Figures \ref{fig:dta_1} and \ref{fig:post_Km_11} we note that the model with $\kappa = 10$ ($\msf{B}$) mixes better the components of the mixture than the other ones in the sense that fewer components were needed in order to capture every mode. Overall, the cases $\kappa = 10$ ($\msf{B}$) and $\kappa = 100$ ($\msf{C}$), seem to inherit desirable properties from the limiting cases, i.e. $\kappa=0$ ($\msf{A}$) and $\kappa=\infty$ ($\msf{D}$). From the Dirichlet process they  inherit a more efficient component mixing, while from the Geometric process they inherit the flexibility to adapt even if the parameter $\theta$ is not carefully chosen.

\subsection{Analysis for BBSB mixtures with random dependence parameter}

The main objective of this analysis is to determine the optimal value of $\kappa$ for different datasets. To this aim, we first we consider a very simple data set $(\msf{database}$ $\msf{2})$ consisting of $200$ observations that were sampled from a mixture of two Gaussian distributions. And secondly, we examine a more complicated set of data $(\msf{database}$ $\msf{3})$ that contains $200$ observations sampled from a mixture of seven Gaussian kernels with distinct means, variances and weights, this database was created and studied before by  \cite{LMP07}. For each, $\msf{database}$ $\msf{2}$ and $\msf{database}$ $\msf{3}$, we study three BBSB mixtures with parameters $\alpha$ and $\theta$ fixed to distinct values, and compare the estimations with the ones provided by a Pitman-Yor mixture. Recall that this two-parameter generalization of the Dirichlet process has stick-breaking representation with independent l.v.'s $\v_i \sim \msf{Be}(1-\sigma,\theta+i\sigma)$ where $0 \leq \sigma < 1$ and $\theta > -\sigma$ \citep[see for instance][for further details]{PPY92,PY92,P06}. In particular the Dirichlet process is recovered when $\sigma = 0$. For this mixture we fixed $\theta$ and consider the other parameter random with a uniform distribution over $[0,1]$, this way the model is allowed to choose the best value of $\sigma$ for the data set. In a similar spirit, for every BBSB mixture considered here, the parameter $\bkap$ was considered random with a uniform prior distribution over $\{0,1,\ldots,100\}$.

\subsubsection{Results for $\msf{database}$ $\msf{2}$}

\begin{figure}[H]
\centering
\includegraphics[scale=0.39]{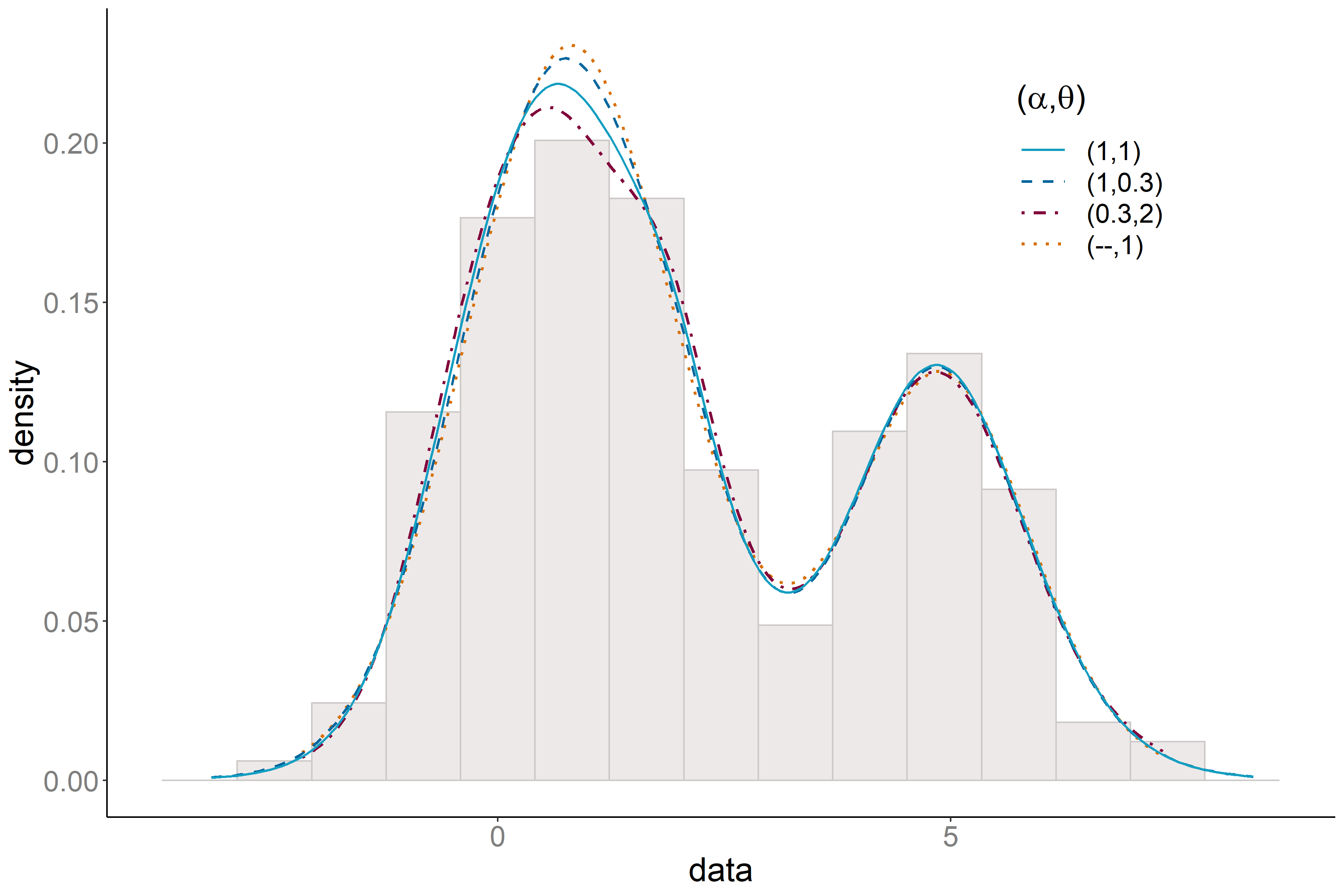}
\caption{\label{fig:dta_2} Estimated densities for $\msf{database}$ $\msf{2}$, taking into account $5000$ iterations of the Gibbs sampler after a burn-in period of $3000$, for three distinct BBSB mixtures with parameters $(\alpha,\theta)$ fixed to $(1,1),(1,0.3)$ and $(0.3,2)$, and a Pitman-Yor mixture with parameter $\theta = 1$.}
\end{figure}

\begin{figure}[H]
\centering
\includegraphics[scale=0.39]{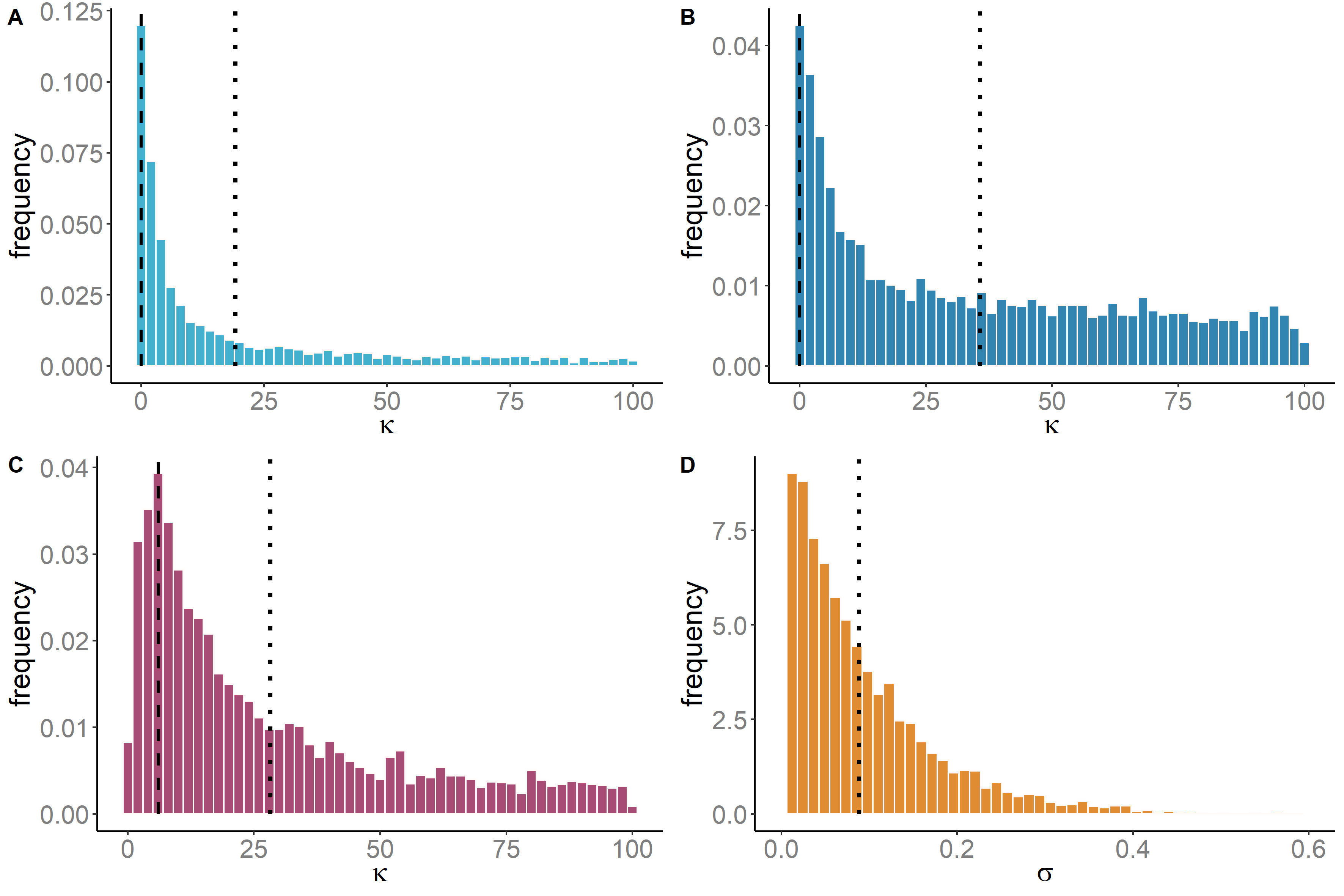}
\caption{\label{fig:kappa_2} Posterior distributions of $\bkap$ $(\msf{A},\msf{B}$ and $\msf{C})$ for the BBSB mixtures with parameters $(\alpha,\theta)$ fixed to $(1,1),(1,0.3)$ and $(0.3,2)$, respectively.
$\msf{D}$ illustrates the posterior distribution of $\bsig$ for the Pitman-Yor mixture with $\theta = 1$. The dotted and dashed lines indicate the posterior means and modes, respectively.}
\end{figure}

In Figure \ref{fig:dta_2} we observe that the estimated densities for the four mixtures adjust well to the data and do not differ significantly. In Figure \ref{fig:kappa_2} we see that every posterior distribution is asymmetrical, hence we will estimate the corresponding randomized parameter by the mode rather that the mean. For the BBSB models with parameter $\alpha = 1$ ($\msf{A}$ and $\msf{B}$), the posterior mode of $\bkap$ equals $0$, suggesting that for this simple data set, the Dirichlet process is an excellent choice. In $\msf{D}$ we see that for  the Pitman-Yor mixture the posterior distribution of $\bsig$ also assigns a bigger probability to values closer to $0$, so it agrees with our models that the Dirichlet process adjust well to this data set. As for the BBSB model with $\alpha = 0.3$ and $\theta = 2$, we observe that the posterior distribution of $\bkap$ ($\msf{C}$) prefers a value bigger than $0$. Explicitly, the posterior mode of this distribution is $\kappa = 6$. {This could be due to the fact that for $\alpha = 0.3$ and $\theta = 2$ the stick-breaking mixture with completely independent l.v.'s is not a good choice for this dataset, so the BBSB mixture corrects this by adjusting the value of the dependence parameter.}

\subsubsection{Results for $\msf{database}$ $\msf{3}$}

\begin{figure}[H]
\centering
\includegraphics[scale=0.35]{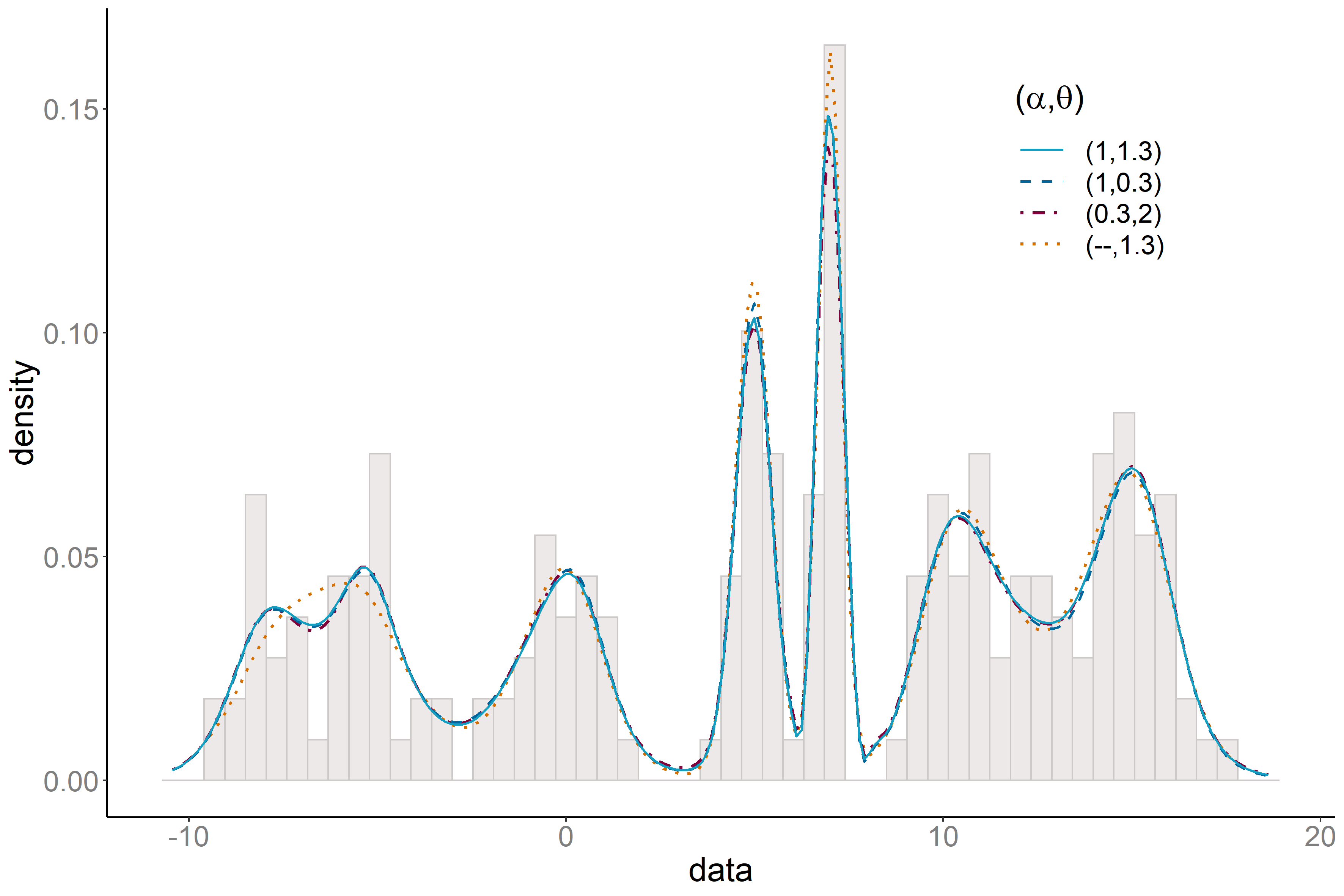}
\caption{\label{fig:dta_3} Estimated densities for $\msf{database}$ $\msf{3}$, taking into account $5000$ iterations of the Gibbs sampler after a burn-in period of $3000$, for three distinct BBSB mixtures with parameters $(\alpha,\theta)$ fixed to $(1,1.3),(1,0.3)$ and $(0.3,2)$, and a Pitman-Yor mixture with parameter $\theta = 1.3$.}
\end{figure}

\begin{figure}[H]
\centering
\includegraphics[scale=0.35]{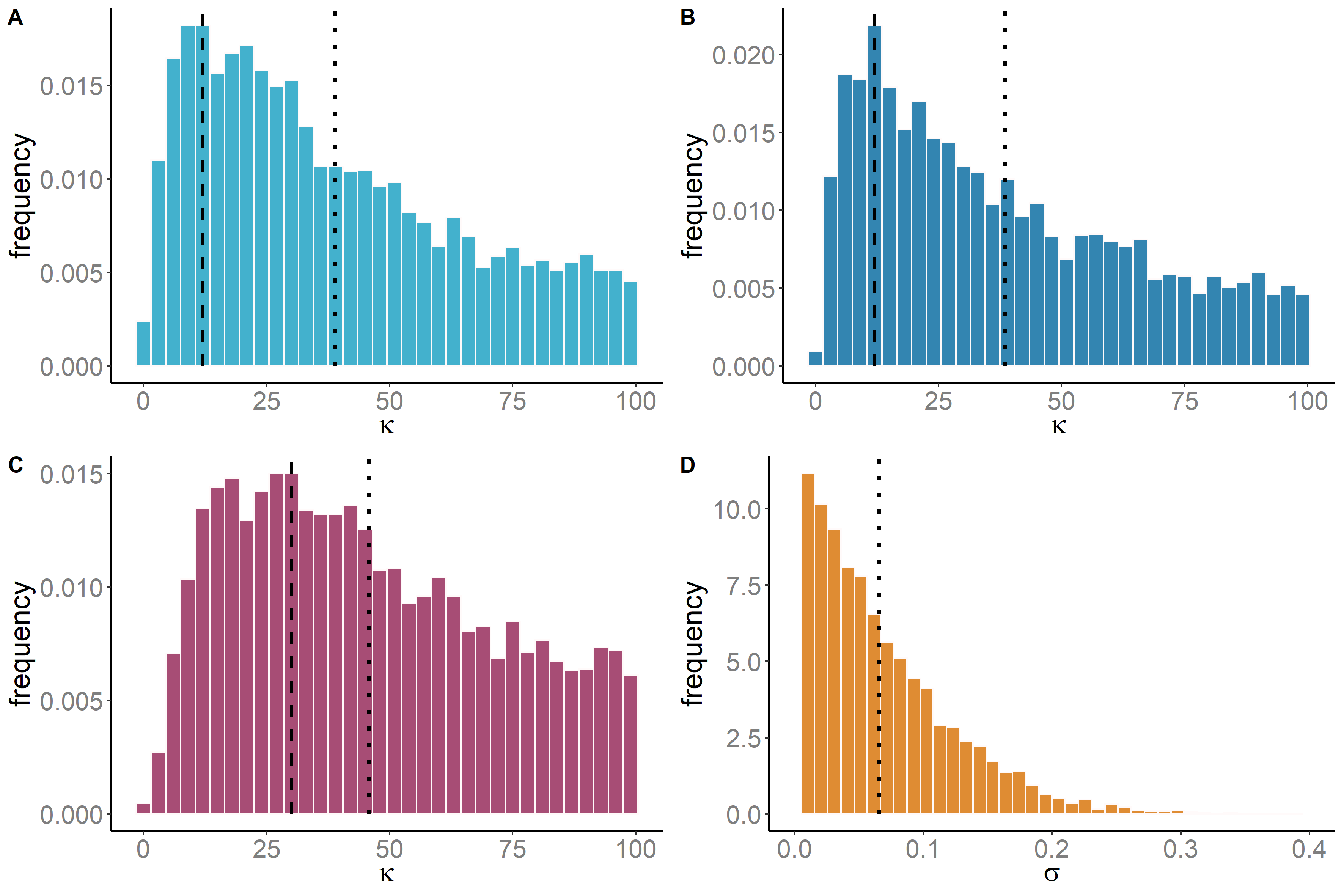}
\caption{\label{fig:kappa_3} Posterior distributions of $\bkap$ $(\msf{A},\msf{B}$ and $\msf{C})$ for the BBSB mixtures with parameters $(\alpha,\theta)$ fixed to $(1,1.3),(1,0.3)$ and $(0.3,2)$, respectively.
$\msf{D}$ illustrates the posterior distribution of $\bsig$ for the Pitman-Yor mixture with $\theta = 1.3$. The dotted and dashed lines indicate the posterior means and modes, respectively.}
\end{figure}

Insomuch as the distributions in Figure \ref{fig:kappa_3} are asymmetrical, once again we estimate the randomized parameter by the posterior modes. In the same figure we observe that the posterior distribution of $\bkap$ for every BBSB mixture $(\msf{A}$, $\msf{B}$ and $\msf{C})$ favours values of $\kappa$ that are bigger than $0$, yet smaller than $50$. Specifically, the posterior modes of $\bkap$ for the BBSB models with $(\alpha,\theta$) fixed to $(1,1.3),(1,0.3)$ and $(0.3,2)$ are $12$, $12$ and $30$, respectively. That is to say, in every case the model estimates that corresponding l.v.'s are dependent. In fact, if we insert the parameters $\alpha = 1,1,0.3$, $\theta = 1.3,0.3,2$, and the posterior mode of $\bkap = 12,12,30$, into Proposition \ref{prop:prop} (d), we estimate the correlation coefficients of consecutive l.v.'s by $0.8992$, $0.9023$ and $0.9288$, respectively. Notice that although the posterior modes of $\bkap$ are not large, these choices affect greatly the dependence of the l.v.'s in question. In particular, for the couple of BBSB mixtures with $\alpha = 1$, this suggest the Dirichlet mixture is not the best choice. Among these two, for the one with $\theta = 1.3$, we chose this parameter so that for the Dirichlet prior $\Esp[\K_{200}] \approx 7$, which coincides with the number of actual modes featured in $\msf{database}$ $\msf{3}$. Even in this case, the posterior distribution of $\bkap$ suggest that other BBSB models fit better than the Dirichlet mixture. As to the Pitman-Yor mixture, for which $\theta$ was also chosen as above, we see in Figure \ref{fig:kappa_3} $(\msf{D})$ that the posterior distribution of $\bsig$ favours values close to $0$. Meaning that this model suggests that among the possibilities, the Dirichlet process is the best fit. However, if we concentrate in Figure \ref{fig:dta_3} we see that the estimated densities by all three BBSB mixtures adjust well the data and recover the seven modes featuring the data set, whilst the Pitman-Yor model confuses the couple of modes in the left hand side of the figure. This suggests the class of BBSB mixtures offers a bigger capacity to adjust to the data by tuning the parameter $\kappa$, than the class of Pitman-Yor mixtures have by tuning the parameter $\sigma$.

\section{Discussion}

By using Beta chains as the l.v.'s of stick-breaking sequences, we were able to construct a new family of distributions over the infinite dimensional simplex, hence a new class of species sampling priors. The parameter, $\kappa$, that modulates the dependence among the elements of the Beta chain, also modulates the ordering of the corresponding weights. While the choice $\kappa = 0$ and $\alpha = 1$ recovers the size-biased permutation of the weights of Dirichlet processes, as $\kappa \to \infty$, we recover the decreasing ordered weights of Geometric processes, both classes of processes being models of interest. This approach to define priors also allows the construction of random measures that are hybrids between Dirichlet and Geometric processes. Furthermore, how similar is the BBSB prior to one model or the other can  also be tuned by the parameter $\kappa$. As to the prior distribution of $\K_n$, generally speaking, we found that a larger value of $\kappa$ translates to a less informative prior. This in turn allows more flexible models in a density estimation context. In general the class of BBSB mixtures offers models with a great capacity to adapt to distinct data sets and models with a efficient component mixing. By endowing the parameter $\kappa$ with a prior distribution, one can estimate its optimal value for a given data set, thus choose the BBSB mixture that admits the optimal balance between flexibility and efficient mixing.

The present work gives rise to interesting questions, such as how to characterize the distribution of $\K_n$ for BBSB priors and analyse its asymptotic behaviour as $n \to \infty$, or even further study the underlying exchangeable partition probability functions. As to the orderings of the weights, it is also of interest to compute or approximate $\Prob[\w_j > \w_{j+1}]$ for a fixed value of $\kappa$, and to determine the rate at which $\Prob[\w_j > \w_{j+1}] \to 1$ as $\kappa \to \infty$. On a non-exchangeable context \citep[e.g.][]{GL17,DMRM04}, one could also use the Beta-Binomial transition to model dependence between two of more species sampling processes whose weights enjoy the stick-breaking decomposition. Hopefully, the present paper motivates the study of stick-breaking sequences featuring dependent l.v.'s, that might even lead to other type of priors. 

\section{Acknowledgements} The first author gratefully thanks the support of CONACyT PhD scholarship program and CONACyT project 241195. The second author gratefully acknowledges the support of CONTEX project 2018-9B as well as the hospitality of the University of Bath, where part of the project was done, during a {\it Global Professor} research visit.

\appendix

\gdef\thesection{\Alph{section}}
\makeatletter
\renewcommand\@seccntformat[1]{Appendix \csname the#1\endcsname.\hspace{0.5em}}
\makeatother

\section{}

\subsection{Convergence of probability measures}
To formally give the proof of the main results, we recall some topological details of measure spaces. For a Polish space $S$, with Borel $\sigma$-algebra $\B(S)$, we denote by $\cl{P}(S)$ the space of all probability measures over $(S,\B(S))$. A well-known metric on $\cl{P}(S)$ is the L\'evy-Prokhorov metric given by
\begin{equation}\label{eq:weak_metric}
d_L(P,P') = \inf\{\varepsilon > 0: P(A) \leq P'\l(A^{\varepsilon}\r) + \varepsilon, P'(A) \leq P\l(A^{\varepsilon}\r)+ \varepsilon, \A A \in \B(S)\},
\end{equation}
for any $P,P' \in \cl{P}(S)$, and where $A^{\varepsilon} = \{s \in S: d(s,A) < \varepsilon\}$, $d(s,A) = \inf\{d(a,s): a \in A\}$ and $d$ is some complete metric on $S$. For probability measures $P,P_1,P_2,\ldots$ it is said that $P_n$ converges weakly to $P$, denoted by $P_n \wto P$, whenever $\int_S f dP_n \to \int_S f dP$ for every continuous bounded function $f:S \to [0,\infty)$. This condition is known to be equivalent to $d_L(P_n,P) \to 0$, and to $\bgamma_n \dto \bgamma$, whenever $\bgamma_n \sim P_n$ and $\bgamma \sim P$. $\cl{P}(S)$, equipped with the topology of weak convergence, is Polish again. Its Borel $\sigma$-field, $\B(\cl{P}(S))$, can equivalently be defined as the $\sigma$-algebra generated by all the projection maps $\{P \mapsto P(B): B \in \B(S)\}$. In this sense the random probability measures (measurable mappings from a probability space $(\Omega,\F,\Prob)$ into $(\cl{P}(S), \B(S))$), $\bmu,\bmu_1,\bmu_2,\ldots$, are said to converge weakly, a.s.  whenever $\bmu_n(\omega) \wto \bmu(\omega)$, for every $\omega$ outside a $\Prob$-null set. Analogously, if $\int_S f d\bmu_n \dto \int_S f d\bmu$ for every continuous bounded function $f:S \to [0,\infty)$, it is said that $\bmu_n$ converges weakly in distribution to $\bmu$, denoted by $\bmu_n \dwto \bmu$. Evidently, $\bmu_n \wto \bmu$ a.s. implies $\bmu_n \dwto \bmu$, which, in turns is a necessary and sufficient condition for $\bmu_n \dto \bmu$. For further details see for instance \cite{P67}, \cite{B68} or \cite{K17}.

\subsection{Proof of Proposition \ref{prop:prop}}\label{ap:prop:prop}

a) Using elementary properties of conditional expectation and the fact that given $\x_i$, $\v_{i+1}$ is conditionally independent of $\v_i$, we obtain
\[
\Esp[\v_{i+1}|\v_i] = \Esp[\Esp[\v_{i+1}|\x_i]|\v_i] = \Esp\left[\frac{\alpha+\x_i}{\alpha+\theta+\kappa}\bigg|\v_i\right] = \frac{\alpha+\kappa\v_i}{\alpha+\theta+\kappa}.
\]
b) Notice that 
%\[
$
\msf{Var}(\v_{i+1}|\v_i) = \Esp[\msf{Var}(\v_{i+1}|\x_i)|\v_i] + \msf{Var}(\Esp[\v_{i+1}|\x_i]|\v_i)$,
%\]
%we first compute
with
\begin{align*}
\msf{Var}(\Esp[\v_{i+1}|\x_i]|\v_i) = \msf{Var}\left(\frac{\alpha+\x_i}{\alpha+\theta+\kappa}\bigg|\v_i\right) = \frac{\v_i(1-\v_i)\kappa}{(\alpha+\theta+\kappa)^2}.
\end{align*}
Now, note that
\begin{align*}
\Esp[(\alpha+\x_i)(\theta+\kappa-\x_i)|\v_i] &= \msf{Cov}(\alpha+\x_i,\theta+\kappa-\x_i|\v_i)\\&\,\,\,\,+ \Esp[\alpha+\x_i|\v_i]\Esp[\theta+\kappa-\x_i|\v_i]\\
& = -\msf{Var}(\x_i|\v_i) + (\alpha+\kappa\v_i)(\theta+\kappa-\kappa\v_i)\\
& = -\kappa\v_i(1-\v_i) + (\alpha+\kappa\v_i)(\theta+\kappa(1-\v_i))
\end{align*}
Hence
\begin{align*}
\Esp[\msf{Var}(\v_{i+1}|\x_i)|\v_i] &= \Esp\left[\frac{(\alpha+\x_i)(\theta+\kappa-\x_i)}{(\alpha+\theta+\kappa)^2(\alpha+\theta+\kappa+1)}\bigg|\v_i\right] \\&= \frac{-\kappa\v_i(1-\v_i) + (\alpha+\kappa\v_i)(\theta+\kappa(1-\v_i))}{(\alpha+\theta+\kappa)^2(\alpha+\theta+\kappa+1)},
\end{align*}
and we can conclude the proof of b),
\begin{align*}
\msf{Var}(\v_{i+1}|\v_i) %& = \frac{-\kappa\v_i(1-\v_i) + (\alpha+\kappa\v_i)(\theta+\kappa(1-\v_i))+\v_i(1-\v_i)\kappa(\alpha+\theta+\kappa+1)}{(\alpha+\theta+\kappa)^2(\alpha+\theta+\kappa+1)}\\
& = \frac{(\alpha+\kappa\v_i)(\theta+\kappa(1-\v_i))+\kappa\v_i(1-\v_i)(\alpha+\theta+\kappa)}{(\alpha+\theta+\kappa)^2(\alpha+\theta+\kappa+1)}.
\end{align*}
c) We first note that as a consequence of the joint reversibility of the Beta-Binomial chain, $\v_i \sim \msf{Be}(\alpha+\x_i,\theta+\kappa-\x_i)$ conditionally given $\x_i$, thus
\begin{align*}
\Esp[\v_i\v_{i+1}] & = \Esp[\Esp[\v_i\v_{i+1}|\x_i]] = \Esp[\Esp[\v_i|\x_i]\Esp[\v_{i+1}|\x_i]] = \Esp\left[\left(\frac{\alpha+\x_i}{\alpha+\theta+\kappa}\right)^2\right],
\end{align*}
conditioning on $\v_i$, we obtain
\begin{align*}
\Esp\left[\left(\frac{\alpha+\x_i}{\alpha+\theta+\kappa}\right)^2\right] & = \Esp\left[\Esp\left[\left(\frac{\alpha+\x_i}{\alpha+\theta+\kappa}\right)^2\bigg|\v_i\right]\right]\\
& =  \Esp\left[\frac{\alpha^2 + 2\alpha\Esp[\x_i|\v_i]+ \Esp[\x_i^2|\v_i]}{(\alpha+\theta+\kappa)^2}\right]\\
& = \frac{\alpha^2 + 2\alpha \kappa\Esp[\v_i]+ \kappa\Esp[\v_i] + \kappa(\kappa-1)\Esp[\v_i^2]}{(\alpha+\theta+\kappa)^2}\\
& = \l[\alpha^2+\frac{\kappa(2\alpha^2+\alpha)}{\alpha+\theta}+\frac{\kappa(\kappa-1)\alpha(\alpha+1)}{(\alpha+\theta)(\alpha+\theta+1)}\r](\alpha+\theta+\kappa)^{-2},
\end{align*}
hence
\begin{align*}
\msf{Cov}(\v_i,\v_{i+1})  = \Esp[\v_i\v_{i+1}]-\Esp[\v_i]\Esp[\v_{i+1}]
%& = (\alpha+\theta+\kappa)^{-2}\l[\alpha^2+\frac{\kappa(2\alpha^2+\alpha)}{\alpha+\theta}+\frac{\kappa(\kappa-1)\alpha(\alpha+1)}{(\alpha+\theta)(\alpha+\theta+1)}\r]-\frac{\alpha^2}{(\alpha+\theta)^2}\\
 = \frac{\kappa\alpha\theta}{(\alpha+\theta)^2(\alpha+\theta+1)(\alpha+\theta+\kappa)}.
\end{align*}
d) The correlation simplifies as follows
\[
\rho_{\v_i,\v_{i+1}} = \frac{\msf{Cov}(\v_i,\v_{i+1})}{\sqrt{\msf{Var}(\v_i)}\sqrt{\msf{Var}(\v_{i+1})}} 
%= \frac{\kappa\alpha\theta(\alpha+\theta)^2(\alpha+\theta+1)}{\alpha\theta(\alpha+\theta)^2(\alpha+\theta+1)(\alpha+\theta+\kappa)} 
= \frac{\kappa}{\alpha+\theta+\kappa}.
\]

\subsection{Proof of Proposition \ref{prop:BB_DP_G}}\label{ap:prop:BB_DP_G}

To prove Proposition~\ref{prop:BB_DP_G} we need some preliminary results.

\begin{lem}[Continuous mappings]\label{lem:cm}
Let $S$ and $T$ be Polish spaces. Let $\bm{\eta},\bm{\eta}_1,\ldots$ be random elements taking values in $S$, with $\bm{\eta}_n \dto \bm{\eta}$, and consider some measurable mappings $f,f_1,f_2\ldots$ from $S$ into $T$ satisfying $f_n(s_n) \to f(s)$, for every $s_n \to s$ in $S$. Then $f_n(\bm{\eta}_n) \dto f(\bm{\eta})$. 
\end{lem}

\begin{lem}\label{lem:joint_conv}
Let $\bgamma^{n} = (\bgamma^n_1, \bgamma^n_2,\ldots)$, $\bgamma = (\bgamma_1,\bgamma_2,\ldots)$ be random sequences taking values in a Polish space $S$. Then $\bgamma^n \dto \bgamma$ if and only if 
\[
(\bgamma^n_1,\ldots,\bgamma^{n}_i) \dto (\bgamma_1,\ldots,\bgamma_i), \text{ for every } i \geq 1.
\]
\end{lem}

Lemmas \ref{lem:joint_conv} and \ref{lem:cm} are well-known result in probability theory, see for instance Theorems 4.27 and 4.29, respectively,  in  \cite{K02}.

\begin{lem}\label{lem:rcd_conv}
Let $S$ and $T$ be Polish spaces. Consider some random elements $\bgamma,\bgamma_1,\bgamma_2,\ldots$ and $\bm{\eta},\bm{\eta}_1,\bm{\eta}_2,\ldots$ taking values in $S$ and $T$, respectively. Let $\rho$ be the distribution of $\bgamma$ and $\rho_n$ the distribution of $\bgamma_n$, also consider some regular versions, $\pi(\cdot|\bgamma)$ and $\pi_n(\cdot|\bgamma_n)$, of $\Prob[\bm{\eta} \in \cdot\,|\bgamma]$ and $\Prob[\bm{\eta}_n \in \cdot\,|\bgamma_n]$ respectively. If $\rho_n \wto \rho$ and for every $s_n \to s$ in $S$ we have that $\pi_n(\cdot|s_n) \wto \pi(\cdot|s)$, then $(\bgamma_n,\bm{\eta}_n) \dto (\bgamma,\bm{\eta})$.
\end{lem}

\textbf{Proof:}

Let $g:S\times T \to \R$ be a continuous and bounded function.
Define $f,f_1,f_2,\ldots:S \to \R$ by
\begin{align*}
f_n(s) = \int g(s,t)\pi_n(dt|s) \quad \text{ and } \quad f(s) = \int g(s,t)\pi(dt|s)
\end{align*}
The first thing we will prove is that 
\begin{equation}\label{eq:fnf}
f_n(s_n) \to f(s) \quad \text{ as } \quad s_n \to s.
\end{equation}
So let $s_n \to s$. Choose some random elements $\bm{\zeta},\bm{\zeta}_1,\bm{\zeta}_2,\ldots$ with $\bm{\zeta}_n \sim \pi_n(\cdot|s_n)$ and $\bm{\zeta} \sim \pi(\cdot|s)$, this way, $\bm{\zeta}_n \dto \bm{\zeta}$ by hypothesis. Define $h,h_1,h_2,\ldots:T \to \R$ by $h_n(t) = g(s_n,t)$ and $h(t) = g(s,t)$. As $g$ is continuous, we have that $h_n(t_n) = g(s_n,t_n) \to g(s,t) = h(t)$, for every $t_n \to t$ in $T$. By Lemma \ref{lem:cm} we obtain $h_n(\bm{\zeta}_n) \dto h(\bm{\zeta})$, which in turn implies
\begin{align*}
\int g(s_n,t)\pi_n(dt|s_n) &= \Esp[g(s_n,\bm{\zeta}_n)] = \Esp[h_n(\bm{\zeta}_n)]\\ &\to \Esp[h(\bm{\zeta})] = \Esp[g(s,\bm{\zeta})] = \int g(s,t)\pi(dt|s).
\end{align*}
Since $s_n \to s$ was arbitrary, this proves equation \eqref{eq:fnf}, which together with the hypothesis and by Lemma \ref{lem:cm} show that $f_n(\bgamma_n) \dto f(\bgamma)$. Particularly,
\begin{align}\label{eq:EfnEf}
\int\l(\int g(s,t)\pi_n(dt|s)\r)&\rho_n(ds) = \Esp[f_n(\bgamma_n)] \nonumber\\
&\to \Esp[f(\bgamma)] = \int\l(\int g(s,t)\pi(dt|s)\r)\rho(ds).
\end{align}
Note that the double integral in the left side of equation \eqref{eq:EfnEf} coincides with $\Esp[g(\bgamma_n,\bm{\eta_n})]$, whilst the one at the right side coincides with $\Esp[g(\bgamma,\bm{\eta})]$. That is, we have proven that $\Esp[g(\bgamma_n,\bm{\eta_n})] \to \Esp[g(\bgamma,\bm{\eta})]$, for every continuous and bounded function $g:S\times T \to \R$. Or equivalently $(\bgamma_n,\bm{\eta}_n) \dto (\bgamma,\bm{\eta})$.

\begin{lem}\label{lem:Bin_L2}
Let $(\x_n)_{n \geq 1}$ be a sequence of random variables such that $\x_n \sim \msf{Bin}(n,p_n)$ for every $n \geq 1$ and where $p_{n} \to p$ in $[0,1]$. Then 
\[
\frac{\x_n}{n} \quad \overset{\cl{L}_2}{\to} \quad p.
\]
\end{lem}

\textbf{Proof:}

For $n \geq 1$, 
\begin{equation}\label{eq:bin_L2}
\begin{aligned}
\Esp\l[\l(\frac{\x_n}{n}-p\r)^2\r] & = \frac{1}{n^2}\Esp\l[\x_n^2\r] -\frac{2p}{n}\Esp[\x_n] + p^2\\
%& = \frac{{n}p_{n}(1-p_{n})+n^2p_{n}^2}{n^2} - 2pp_{n} + p^2\\
& = \frac{p_{n}(1-p_{n})}{n} + (p_{n}-p)^2.
\end{aligned}
\end{equation}
By taking limits as $n \to \infty$ in \eqref{eq:bin_L2} we obtain
\begin{equation*}
\lim_{n \to \infty} \Esp\l[\l(\frac{\x_n}{n}-p\r)^2\r] = 0.
\end{equation*}

\textbf{Proof of Proposition \ref{prop:BB_DP_G}:}

(i) Insomuch as the corresponding spaces are Borel, we may construct on some probability space $(\hat{\Omega},\hat{\F},\hat{\Prob})$ a Beta-Binomial chain $(\hat{\V},\hat{\X})$ with parameters $(0,\alpha,\theta)$. Now, the elements of $\hat{\V}$ are conditionally independent given $\hat{\X}$, and given that $\kappa = 0$, $\hat{\X} \aseq (0,0,\ldots)$, so we may think of $\hat{\X}$ as if it was deterministic, which implies that the elements of $\hat{\V}$ must be independent and $\msf{Be}(\alpha,\theta)$ distributed.

(ii): For every $\kappa \geq 1$, let $\V^{(\kappa)} = \l(\v^{(\kappa)}_i\r)_{i \geq 1}$ be a Beta chain with parameters $(\kappa,\alpha,\theta)$, and let $\pi_{\kappa}\l(\cdot\big|\v^{(\kappa)}_i\r)$ be some regular version of $\Prob\l[\v^{(\kappa)}_{i+1} \in \cdot\,\big|\v^{(\kappa)}_i\r]$ (which clearly does not depends on $i$). Further let $\bm{\lambda}\sim \msf{Be}(\alpha,\theta)$ and fix $\pi(\cdot|\bm{\lambda}) = \delta_{\bm{\lambda}}$. The first thing we are interested in proving is that for every $p_\kappa \to p$ in $[0,1]$ we have that 
\begin{equation}\label{eq:rcd_conv}
\pi_\kappa(\cdot|p_\kappa) \wto \pi(\cdot|p).
\end{equation}
So, let $p_\kappa \to p$ in $[0,1]$, by Lemma \ref{lem:Bin_L2} and given that all the corresponding spaces are Borel, we may construct on a probability space $(\hat{\Omega},\hat{\F},\hat{\Prob})$, with expectations $\hat{\Esp}[\cdot]$, some pairs of r.v.'s $\l(\hat{\x}_\kappa,\hat{\v}_\kappa\r)_{\kappa \geq 1}$ such that $\hat{\x}_\kappa \sim \msf{Bin}(\kappa,p_\kappa)$, $\l\{\hat{\v}_\kappa|\hat{\x}_\kappa\r\} \sim \msf{Be}(\alpha+\hat{\x}_\kappa, \theta+\kappa-\hat{\x}_\kappa)$, and $\hat{\x}_\kappa/\kappa$ $\asto$ $p$. Note that marginally $\hat{\v}_\kappa \sim \pi_\kappa(\cdot|p_\kappa)$ so to prove equation \eqref{eq:rcd_conv}, it suffices to show $\hat{\v}_\kappa \dto p$.

Conditionally given $\hat{\x}_\kappa$ the moment generator function of $\hat{\v}_\kappa$ is
\begin{equation}\label{eq:mgf_v|x}
\hat{\Esp}\l[e^{t\hat{\v}_\kappa}\big|\hat{\x}_\kappa\r] = 1 + \sum_{k=1}^{\infty}\l(\prod_{r=0}^{k-1}\frac{\alpha+\hat{\x}_\kappa+r}{\alpha+\theta+\kappa+r}\r)\frac{t^k}{k!}, \quad t \in \R.
\end{equation}
By construction we have that $\hat{\x}_\kappa/\kappa$ $\asto$ $p$, which means that for every $r \geq 0$,
\begin{equation}\label{eq:mgf_v|x_lim}
\frac{\alpha+\hat{\x}_\kappa+r}{\alpha+\theta+\kappa+r} = \l(\frac{\alpha+r}{\kappa}+\frac{\hat{\x}_\kappa}{\kappa}\r)\l(\frac{\alpha+\theta+r}{\kappa}+1\r)^{-1} \asto p,
\end{equation}
as $\kappa \to \infty$. Hence by the tower property of conditional expectation, equations \eqref{eq:mgf_v|x} and \eqref{eq:mgf_v|x_lim}, and Lebesgue dominated convergence theorem (the corresponding functions are dominated by $e^t$) we obtain
\begin{align*}
\lim_{\kappa \to \infty} \hat{\Esp}\l[e^{t\hat{\v}_\kappa}\r] & = \lim_{\kappa \to \infty} \hat{\Esp}\l[\hat{\Esp}\l[e^{t\hat{\v}_\kappa}|\hat{\x}_\kappa\r]\r]\\
& = \hat{\Esp}\l[1 + \sum_{k=1}^{\infty}\l(\prod_{r=0}^{k-1}\lim_{\kappa \to \infty}\frac{\alpha+\hat{\x}_\kappa+r}{\alpha+\theta+\kappa+r}\r)\frac{t^k}{k!}\r]\\
& = \hat{\Esp}\l[1 + \sum_{k=1}^{\infty}\frac{(pt)^k}{k!}\r]\\
& = e^{tp},
\end{align*}
which proves altogether  $\hat{\v}_\kappa \dto p$ and equation \eqref{eq:rcd_conv}.\\

Returning to the original Beta chains, we have that $\v^{(\kappa)}_1 \deq \bm{\lambda}$ for every $\kappa \geq 1$, so trivially, $\v^{(\kappa)}_1 \dto \bm{\lambda}$, this together with equation \eqref{eq:rcd_conv} and the recursive application of Lemma \ref{lem:rcd_conv} allows us to obtain
\[
\l(\v^{(\kappa)}_1,\ldots,\v^{(\kappa)}_i\r) \dto (\bm{\lambda},\ldots,\bm{\lambda}), \quad i \geq 1,
\]
and by Lemma \ref{lem:joint_conv} we can conclude $\V^{(\kappa)} = \l(\v^{(\kappa)}_i\r)_{i \geq 1} \dto (\bm{\lambda},\bm{\lambda},\ldots)$.\\ 

\subsection{Proof of Proposition \ref{prop:BB_s.b.s.}}\label{ap:prop:BB_s.b.s.}

For sequences that enjoy the decomposition \eqref{eq:s.b.s} we may equivalently prove that 
\[
\l(1 - \sum_{i=1}^{j} \w_i\r) = \prod_{i=1}^{j}(1-\v_i) \asto 0,
\]
as $j \to \infty$ \citep[see for instance][]{GvdV17}. Further, these r.v.'s are non-negative and bounded by $1$, thus it is enough to show that 
\begin{equation}\label{eq:esp_prod}
\lim_{j \to \infty}\Esp\l[\prod_{i=1}^{j}(1-\v_i)\r]  = 0.
\end{equation}
As the corresponding spaces are Borel, (after possibly enlarging the original probability space) it is possible to construct a Binomial chain $\X$ such that $(\V,\X)$ defines a Beta-Binomial chain. Conditionally given $\X = \{\x_i\}_{i \geq 1}$, the elements of $\V = \{\v_i\}_{i \geq 1}$ are independent with, $\{\v_1|\x_1\} \sim \msf{Be}(\alpha + \x_1,\theta + \kappa -\x_1)$ and $\l\{\v_{i+1}|\x_i,\x_{i+1}\r\} \sim \msf{Be}(\alpha + \x_i + \x_{i+1},\theta+2\kappa-\x_{i}-\x_{i+1})$, for $i \geq 1$. Hence
\begin{align*}
\Esp\l[\prod_{i=1}^{j}(1-\v_i)\r] & = \Esp\l[\Esp\l[\prod_{i=1}^{j}(1-\v_i)\bigg|\X\r]\r]\\
& = \Esp\l[\Esp[(1-\v_1)|\x_1]\prod_{i=2}^{j}\Esp\l[(1-\v_i)|\x_{i-1},\x_i\r]\r]\\
& = \Esp\l[\frac{\theta+\kappa-\x_1}{\alpha+\theta+\kappa}\prod_{i=2}^{j}\frac{\theta+2\kappa-\x_{i}-\x_{i-1}}{\alpha+\theta+2\kappa}\r].
\end{align*}
Recalling that $0 \leq \x_i \leq \kappa$ a.s. we obtain 
\[
\frac{\theta}{\alpha+\theta+\kappa}\l(\frac{\theta}{\alpha+\theta+2\kappa}\r)^{j-1} \leq \Esp\l[\prod_{i=1}^{j}(1-\v_i)\r] \leq \frac{\theta+\kappa}{\alpha+\theta}\l(\frac{\theta+2\kappa}{\alpha+\theta+2\kappa}\r)^{j-1},
\]
for every $j \geq 1$. Finally by taking limits as $j \to \infty$ in the last equation, \eqref{eq:esp_prod} follows.

\subsection{Proof of Theorem \ref{theo:BB_DP_G}}\label{ap:theo:BB_DP_G}

To prove Theorem \ref{theo:BB_DP_G} we will first prove a couple of elementary results. 

\begin{lem}\label{lem:dL_bound}
Let $S$ be a Polish space and fix some distinct $s_1,s_2,\ldots \in S$, let $p = (p_1,p_2,\ldots)$ and $q = (q_1,q_2,\ldots)$ be elements of $\Delta_{\infty}$ and define $P = \sum_{j \geq 1}p_j\delta_{s_j}$ and $Q = \sum_{j \geq 1}q_j\delta_{s_j}$. Then for $d_L$ as in equation \eqref{eq:weak_metric}
\[
d_L(P,Q) \leq \sum_{j \geq 1}|p_j-q_j|.
\]
\end{lem}

\textbf{Proof:}\\

Define $\varepsilon(p,q) = \sum_{j \geq 1}|p_j-q_j|$, by definition of $d_L$, it suffices to prove for all $A \in \B(S)$
\begin{equation}\label{eq:epsilon_condition}
P(A) \leq Q\l(A^{\varepsilon(p,q)}\r) + \varepsilon(p,q), \quad \text{ and }\quad
Q(A) \leq P\l(A^{\varepsilon(p,q)}\r) + \varepsilon(p,q), %\quad \quad \A A \in \B(S),
\end{equation}
So let $A \in \B(S)$ and set $M_A = \{j \geq 1: s_j \in A\}$, then
\begin{align*}
P(A) = \sum_{j \in M_A}P(\{s_j\}) = \sum_{j \in M_A}p_j & \leq \sum_{j \in M_A}q_j + \sum_{j \in M_A}|p_j-q_j|\\
& \leq Q(A) + \varepsilon(p,q)\\
& \leq Q\l(A^{\varepsilon(p,q)}\r) + \varepsilon(p,q).
\end{align*}
Analogously, we have that $Q(A) \leq P\l(A^{\varepsilon(p,q)}\r) + \varepsilon(p,q)$.

\begin{lem}\label{lem:dL_cont_map}
For fixed and distinct elements $s_1,s_2,\ldots \in S$, the mapping,
\[
(w_1,w_2,\ldots) \mapsto \sum_{j \geq 1}w_j\delta_{s_j},
\]
from $\Delta_{\infty}$ into $\cl{P}(S)$ is continuous with respect to the weak topology.
\end{lem}

\textbf{Proof:}\\

Let $w^{(n)} = \l(w^{(n)}_1,w^{(n)}_2,\ldots\r)$ and $w = (w_1,w_2,\ldots)$ be any elements of $\Delta_{\infty}$ such that $w^{(n)}_j \to w_j$, for every $j \geq 1$. Define $P^{(n)} = \sum_{j \geq 1}w^{(n)}_j\delta_{s_j}$ and $P = \sum_{j \geq 1}w_j\delta_{s_j}$. By Lemma \ref{lem:dL_bound}
\begin{equation*}
d_L\l(P^{(n)},P\r) \leq \sum_{j \geq 1}|w^{(n)}_j-w_j| \leq \sum_{j \geq 1}w^{(n)}_j + \sum_{j \geq 1}w_j = 2,
\end{equation*}
and by the general Lebesgue dominated convergence theorem we obtain
\[
\lim_{n \to \infty}d_L\l(P^{(n)},P\r) = \lim_{n \to \infty}\sum_{j \geq 1}|w^{(n)}_j-w_j| = \sum_{j \geq 1}\lim_{n \to \infty}|w^{(n)}_j-w_j| = 0,
\]
which means that the mapping $(w_1,w_2,\ldots) \mapsto \sum_{j \geq 1}w_j\delta_{s_j}$ is continuous.

\begin{rem}
Despite the choice of the metric, $\rho$, in $\Delta_{\infty}$, as long as $\rho$ generates the Borel $\sigma$-algebra, $\rho\l(w^{(n)},w\r) \to 0$ implies $|w^{(n)}_j-w_j| \to 0$, for every $j \geq 1$. For this reason, in the above proof we did not discuss the details on the metric, of $\Delta_\infty$, that is being used.
\end{rem}

\textbf{Proof of Theorem \ref{theo:BB_DP_G}:}\\

The proof of (i) follows directly from Proposition \ref{prop:BB_DP_G} (i). To prove (ii), note that by Proposition \ref{prop:BB_DP_G} (ii) and given that all the corresponding spaces are Borel, we can construct on a probability space $(\hat{\Omega}, \hat{\F}, \hat{\Prob})$, Beta chains $\hat{\V}^{(\kappa)} = \l(\hat{\v}^{(\kappa)}_i\r)_{i \geq 1}$ with parameters $(\kappa,\alpha,\theta)$ and a $\hat{\bm{\lambda}} \sim \msf{Be}(\alpha,\theta)$ such that $\hat{\v}^{(\kappa)}_i \asto \hat{\bm{\lambda}}$, as $\kappa \to \infty$, for every $i \geq 1$. Define also an independent sequence, $\hat{\bXi} = \l(\hat{\bxi}_j\r)_{j \geq 1}$, with $\hat{\bxi}_j \iid P_0$. Now, for $\kappa \geq 1$ set
\[
\hat{\w}_j^{(\kappa)} = \hat{\v}^{(\kappa)}_j\prod_{i=1}^{j-1}\l(1-\hat{\v}^{(\kappa)}_i\r), \quad j \geq 1, \quad \text{ and }  \quad \hat{\bmu}^{(\kappa)} = \sum_{j \geq 1}\hat{\w}^{(\kappa)}_j\delta_{\hat{\bxi}_j},
\]
with the empty product equating to $1$, also set $\hat{\bmu} = \sum_{j \geq 1}\hat{\bm{\lambda}}\l(1-\hat{\bm{\lambda}}\r)^{j-1}\,\delta_{\hat{\bxi}_j}$, so 
\begin{equation}\label{eq:mu_hatmu}
\hat{\bmu}^{(\kappa)} \deq \bmu^{(\kappa)}, \quad \kappa \geq 1 \quad \text{ and } \quad \hat{\bmu} \deq \bmu.
\end{equation}
As the mapping $\l(\hat{\v}^{(\kappa)}_1,\ldots,\hat{\v}^{(\kappa)}_j\r) \mapsto \hat{\w}_j^{(\kappa)}$ is continuous, we have that
\[
\hat{\w}_j^{(\kappa)} \asto \hat{\bm{\lambda}}\l(1-\hat{\bm{\lambda}}\r)^{j-1}, \quad j \geq 1.
\]
For the sequence $\hat{\bXi}$, the diffuseness of $P_0$ implies that for $i \neq j$, $\hat{\bxi}_i \neq \hat{\bxi}_j$ a.s., since we are dealing with a countable number of random variables, there exist some $B \in \hat{F}$ such that $\hat{\Prob}[B] = 1$ and for every $\omega \in B$
\[
\hat{\w}_j^{(\kappa)}(\omega) \to\hat{\bm{\lambda}}(\omega)\l(1-\hat{\bm{\lambda}}(\omega)\r)^{j-1}, \quad j \geq 1, \quad \text{ and } \quad \hat{\bxi}_j(\omega) \neq \hat{\bxi}_i(\omega), \quad i \neq j
\]
By Lemma \ref{lem:dL_cont_map} 
\[
\sum_{j \geq 1}\hat{\w}_j^{(\kappa)}(\omega)\delta_{\hat{\bxi}_j(\omega)} \wto \sum_{j \geq 1}\hat{\bm{\lambda}}(\omega)\l(1-\hat{\bm{\lambda}}(\omega)\r)^{j-1}\,\delta_{\hat{\bxi}_j(\omega)}, \quad \omega \in B
\]
that is, $\hat{\bmu}^{(\kappa)} \wto \hat{\bmu}$ a.s., implying $\hat{\bmu}^{(\kappa)} \dto \hat{\bmu}$. Finally, by equation \eqref{eq:mu_hatmu}, the result follows.

\subsection{Proof of Corollary \ref{cor:dec_prob_weights}}\label{ap:cor:dec_prop_weights}

The proof of (i) can be found in Theorem 1 by \cite{P96b}. To prove (ii) note that we may write
\[
\w_1^{(k)} = \v_1^{(\kappa)}, \quad \w_{j+1}^{(k)} = \frac{\v_{j+1}^{(\kappa)}\l(1-\v_j^{(\kappa)}\r)}{\v_j^{(\kappa)}}\w_j^{(\kappa)}, \quad j \geq 1,
\]
hence
\[
\Prob\l[\w_{j+1}^{(\kappa)}<\w_j^{(\kappa)}\r] = \Prob\l[\v_{j+1}^{(\kappa)}\l(1-\v_j^{(\kappa)}\r)<\v_j^{(\kappa)}\r].
\]
By the second part of Proposition \ref{prop:BB_DP_G} and as the corresponding spaces are Borel, we may construct on some probability space, $(\hat{\Omega}, \hat{\F}, \hat{\Prob})$, with expectations $\hat{\Esp}[\cdot]$, Beta chains, $\l(\hat{\v}_i^{(\kappa)}\r)_{i \geq 1}$, with parameters $(\kappa,\alpha,\theta)$, and a $\hat{\bm{\lambda}} \sim \msf{Be}(\alpha,\theta)$ satisfying
\[
\l(\hat{\v}_i^{(\kappa)}\r)_{i \geq 1} \to (\hat{\bm{\lambda}},\hat{\bm{\lambda}},\ldots) \quad \text{ a.s.}
\]
Then for $j \geq 1$, there exist $A \in \hat{\F}$ with $\hat{\Prob}[A] =1$ and such that for every $\omega \in A$, $\hat{\v}_j^{(\kappa)}(\omega) \to \hat{\bm{\lambda}}(\omega)$ and $\hat{\v}_{j+1}^{(\kappa)}(\omega) \to \hat{\bm{\lambda}}(\omega)$. Fix $\omega \in A$, since $\hat{\bm{\lambda}}(\omega)(1-\hat{\bm{\lambda}}(\omega)) < \hat{\bm{\lambda}}(\omega)$, we may choose $\kappa'$ such that for every $\kappa > \kappa'$, $\hat{\v}_{j+1}^{(\kappa)}(\omega)\l(1-\hat{\v}_j^{(\kappa)}(\omega)\r)<\hat{\v}_j^{(\kappa)}(\omega)$. As $\omega$ was chosen arbitrarily in $A$ we have that $\Ind\l\{\hat{\v}_{j+1}^{(\kappa)}\l(1-\hat{\v}_j^{(\kappa)}\r)<\hat{\v}_j^{(\kappa)}\r\} \to 1$ a.s., as $\kappa \to \infty$. Finally, by Lebesgue dominated convergence theorem we obtain
\begin{align*}
\lim_{\kappa \to \infty}\Prob\l[\w_{j+1}^{(\kappa)}<\w_j^{(\kappa)}\r]
& = \lim_{\kappa \to \infty}\Esp\l[\Ind\l\{\v_{j+1}^{(\kappa)}\l(1-\v_j^{(\kappa)}\r)<\v_j^{(\kappa)}\r\}\r]\\ & = \lim_{\kappa \to \infty}\hat{\Esp}\l[\Ind\l\{\hat{\v}_{j+1}^{(\kappa)}\l(1-\hat{\v}_j^{(\kappa)}\r)<\hat{\v}_j^{(\kappa)}\r\}\r]\\
& = \hat{\Esp}\l[\lim_{\kappa \to \infty}\Ind\l\{\hat{\v}_{j+1}^{(\kappa)}\l(1-\hat{\v}_j^{(\kappa)}\r)<\hat{\v}_j^{(\kappa)}\r\}\r]\\
& = 1.
\end{align*}

\bibliographystyle{dcu}
\bibliography{references}

\end{document}